\documentclass[reqno,12pt]{amsart}
\usepackage{a4wide}
\usepackage{amscd}
\usepackage{amsmath}
\usepackage{amsfonts}
\usepackage{amssymb}
\usepackage{multirow}
\usepackage{color}
\usepackage{hyperref}

\newcommand{\ca}{{\mathfrak{a}}}

\newcommand{\ce}{{\mathfrak{e}}}

\newcommand{\cf}{{\mathfrak{f}}}

\newcommand{\cg}{{\mathfrak{g}}}

\newcommand{\ch}{{\mathfrak{h}}}

\newcommand{\ck}{{\mathfrak{k}}}

\newcommand{\co}{{\mathfrak{o}}}

\newcommand{\cp}{{\mathfrak{p}}}

\newcommand{\cq}{{\mathfrak{q}}}

\newcommand{\crr}{{\mathfrak{r}}}

\newcommand{\cs}{{\mathfrak{s}}}

\newcommand{\cu}{{\mathfrak{u}}}

\newcommand{\cz}{{\mathfrak{z}}}

\newcommand{\RR}{\mathbb{R}}

\newcommand{\CC}{\mathbb{C}}

\newcommand{\HH}{\mathbb{H}}
\newcommand{\OO}{\mathbb{O}}

\newcommand{\NN}{\mathbb{N}}

\DeclareMathOperator{\ad}{ad}
\DeclareMathOperator{\Ad}{Ad}

\DeclareMathOperator{\codim}{codim}
\DeclareMathOperator{\Exp}{Exp}
\DeclareMathOperator{\Fix}{Fix}
\DeclareMathOperator{\id}{id}
\DeclareMathOperator{\rk}{rk}

\newtheorem{thm}{Theorem}[section]

\newtheorem{prop}[thm]{Proposition}
\newtheorem{cor}[thm]{Corollary}
\newtheorem{lm}[thm]{Lemma}
\newtheorem{re}[thm]{Remark}

\numberwithin{equation}{section}

\begin {document}

\title{The index conjecture for symmetric spaces}

\author{J\"{u}rgen Berndt}
\address{King's College London, Department of Mathematics, London WC2R 2LS, United Kingdom}
\email{jurgen.berndt@kcl.ac.uk}
\thanks{}

\author{Carlos Olmos}
\address{Facultad de Matem\'atica, Astronom\'ia y F\'isica, Universidad Nacional de C\'ordoba, 
Ciudad Universitaria, 5000 C\'ordoba, Argentina}
\email{olmos@famaf.unc.edu.ar}
\thanks{The second author acknowledges financial support from Famaf, UNC and Ciem, CONICET}


\begin{abstract}
In 1980, Oni\v{s}\v{c}ik (\cite{O80}) introduced the index of a Riemannian symmetric space as the minimal codimension of a (proper) totally geodesic submanifold. He calculated the index for symmetric spaces of rank $\leq 2$, but for higher rank it was unclear how to tackle the problem. In \cite{BO16}, \cite{BO17}, \cite{BO18} and \cite{BOR19} we developed several approaches to this problem, which allowed us to calculate the index for many symmetric spaces. Our systematic approach led to a conjecture, formulated first in \cite{BO16}, for how to calculate the index. The purpose of this paper is to verify the conjecture.
\end{abstract}

\maketitle 

\thispagestyle{empty}

\section {Introduction}

A generic Riemannian manifold does not admit nontrivial totally geodesic submanifolds apart from geodesics (see e.g.\ \cite{M19}). The situation becomes more interesting when considering Riemannian manifolds with many symmetries. A particularly interesting, and important, class of such manifolds are the Riemannian symmetric spaces. The interplay between the geometric theory of Riemannian symmetric spaces and the algebraic theory of semisimple Lie algebras is very fascinating. In our context, the geometric objects of totally geodesic submanifolds correspond to the algebraic objects of Lie triple systems. Unfortunately, the algebraic equations underlying Lie triple systems turn out to be very complicated in general. The series of papers by Klein (\cite{K08},\cite{K09},\cite{K10a},\cite{K10b}) illustrates very well the complicated nature of classifying Lie triple systems for symmetric spaces of rank $2$. Previously, Wolf (\cite{Wo63}) classified totally geodesic submanifolds in symmetric spaces of rank $1$ by geometric methods. In Riemannian symmetric spaces of rank $\geq 3$ some ``standard'' examples of totally geodesic submanifolds are known, but a classification is out of reach with known methods.

In this context, Oni\v{s}\v{c}ik (\cite{O80}) introduced the notion of index for Riemannian symmetric spaces. He defined the \textit{index} $i(M)$ of a Riemannian symmetric space $M$ as the minimal possible codimension of a nontrivial totally geodesic submanifold. The index provides an obstruction for the existence of totally geodesic immersions between symmetric spaces (see Remark \ref{obstruction}). The basic questions are:

\smallskip
\textit{What is the index of a Riemannian symmetric space, and how to determine it?}

\smallskip
In \cite{O80}, Oni\v{s}\v{c}ik used standard Lie algebraic methods to determine the index of irreducible Riemannian symmetric spaces of rank $2$. In our previous work (\cite{BO16},\cite{BO17},\cite{BO18},\cite{BOR19}) we developed several systematic approaches to these two questions. For many symmetric spaces we were able to determine the index, but more importantly a conjecture, the so-called \textit{Index Conjecture for Symmetric Spaces}, emerged. The conjecture was first formulated in \cite{BO16} and relates the index to fixed point sets of involutions on symmetric spaces. The involutions on irreducible Riemannian symmetric spaces $M$ of compact type and their fixed point sets were determined and studied thoroughly by Nagano and Tanaka (\cite{N88}, \cite{N92}, \cite{NT95}, \cite{NT99}, \cite{NT00}). Every connected component $\Sigma$ of such a fixed point set is a totally geodesic submanifold of $M$. Geometrically, the involution is the geodesic reflection of $M$ in $\Sigma$. This is why such a totally geodesic submanifold, given by an involution, is also called a reflective submanifold. Algebraically, a reflective submanifold corresponds to a Lie triple system for which the orthogonal complement is also a Lie triple system. We call such a Lie triple system a reflective Lie triple system. Reflective submanifolds of irreducible Riemannian symmetric spaces of compact type were studied by Leung  in \cite{L75} and \cite{L79}, where one can also find explicit tables concerning their classification. Since all involutions, or equivalently, reflective submanifolds on irreducible Riemannian symmetric spaces of compact type are explicitly known, it is easy to compute the minimal possible codimension of a reflective submanifold, which we call the {\it reflective index} of $M$ and denote by $i_r(M)$. Our conjecture states:

\smallskip
\textbf{Index Conjecture for Symmetric Spaces.} \textit{For every irreducible Riemannian symmetric space $M$ we have $i(M) = i_r(M)$, unless $M = G_2/SO_4$ or $M = G_2^2/SO_4$.}

\smallskip
The situation for $M = G_2/SO_4$ is quite unique. The special unitary group $SU_3$ is a maximal subgroup of $G_2$ and one of its orbits in $G_2/SO_4$ is totally geodesic and isometric to the $5$-dimensional symmetric space $SU_3/SO_3$, and thus $i(M) \leq 3$. In fact, as Oni\v{s}\v{c}ik proved in \cite{O80}, we have $i(M) = 3$. On the other hand, according to Nagano (\cite{N88}), the only nontrivial involutions on $G_2/SO_4$ are the geodesic symmetries, whose nontrivial fixed point sets are $4$-dimensional and locally isometric to $S^2 \times S^2$. Therefore, the reflective index of $M = G_2/SO_4$ is $i_r(M) = 4$. For $G_2^2/SO_4$ the corresponding statements hold via duality between symmetric spaces of compact type and of noncompact type.

The purpose of this paper is to give an affirmative answer to the Index Conjecture for Symmetric Spaces. With the methods that we developed in our previous work on this topic we could verify the conjecture for some series of classical symmetric spaces, all compact simple Lie groups and all exceptional symmetric spaces. However, none of these methods lead to conclusions for the following three series of classical symmetric spaces:
\begin{itemize}
\item[(i)] $M = SU_{2r+2}/Sp_{r+1}$ for $r \geq 3$. Conjecture: $i(M) = 4r$ for $r \geq 4$ and $i(M) = 11$ for $r = 3$ ($\rk(M) = r$).
\item[(ii)] $M = Sp_{2r+k}/Sp_rSp_{r+k}$ for $k \geq 0$ and  $r \geq \max\{3,k +2\}$. Conjecture: $i(M) = 4r$ ($\rk(M) = r$).
\item[(iii)] $M = SO_{2k+2}/U_{k+1}$ for $k \geq 5$. Conjecture: $i(M) = 2k$ ($\rk(M) = \left[\frac{k+1}{2}\right]$).
\end{itemize}
For each of these three series of classical symmetric spaces we develop a new methodology for calculating the index. None of the three methods can be used to prove the index conjecture for the other series of symmetric spaces. We outline the methods here. By $\Sigma$ we always denote a maximal totally geodesic of $M$. 

(i) {\sc The Lagrangian Grassmannians $M = SU_{2r+2}/Sp_{r+1}$.} If $\Sigma$ is maximal and locally reducible, and if a local de Rham factor of $\Sigma$ has a root system that is not of type $A$, then we show that $\Sigma$ is nonsemisimple (Lemma \ref{SigmaA_r}). From previous work we know that in this case $\Sigma$ is reflective. When the root system of all de Rham factors is of type $A$, we prove some estimates involving dimensions and ranks, which show that $\codim(\Sigma) > i_r(M)$ (Proposition \ref{factorsAAA}). We can thus assume that $\Sigma$ is locally irreducible and its root system is not of type $A$. This leads to three possibilities (Lemma \ref{AAA}): (a) $\Sigma$ is an inner symmetric space; (b) The maximum of the multiplicities of all roots is $\leq 2$; (c) $\Sigma$ is locally isometric to the real Grassmannian $\Sigma = SO_{2s+n}/SO_sSO_{s+n}$ with $s \geq 3$ odd and $n \geq 4$ even. For case (a) we show that every maximal totally geodesic submanifold $\Sigma$ of an outer irreducible symmetric space, where $\Sigma$ an inner symmetric space, is reflective (Proposition \ref{innSigma}). For case (b) we prove an estimate involving the number of reflections in the associated Weyl groups (Proposition \ref{cardinal}), which then leads to $\codim(\Sigma) \geq i_r(M)$ (Proposition \ref{mleq2}). For case (c) we construct an involution on $M$ which has a fixed point set of dimension greater than $\dim(\Sigma)$ (Proposition \ref{grassmannian}).

(ii) {\sc The quaternionic Grassmannians $M = Sp_{2r+k}/Sp_rSp_{r+k}$.} We first show that $i(\Sigma) \leq i(M)$ (Proposition \ref{indexSigmaM}). Using this result we can reduce the problem to the case $k = 0$, that is,  $M = Sp_{2r}/Sp_rSp_r$ (Lemma \ref{spreduction}). We can realize the symplectic group $Sp_r$ as a centrosome in $M = Sp_{2r}/Sp_rSp_r$, and since $i(Sp_r) = 4r-4$, we can conclude that $4r-4 \leq i(M) \leq 4r$ (Lemma \ref{lowbound}). We then derive some inequalities for dimension and rank of the isotropy group of $\Sigma$ and its locally irreducible factors (Lemma \ref{isotropybounds}). Using these inequalities we prove the conjecture for $r \in \{3,4,5\}$ using case-by-case methods (Propositions \ref{isp6sp3sp3}, \ref{isp8sp4sp4} and \ref{isp10sp5sp5}). To simplify these case-by-case calculations we develop a general theory for reducible totally geodesic submanifolds with rank $1$ factors (Section \ref{tgsrk1}), which allows us to dismiss many possibilities. For $r \geq 6$ we then develop an inductive argument (Proposition \ref{isp2rsprspr}), for which we prove an estimate for the codimension of a totally geodesic submanifold in the product of two irreducible Riemannian symmetric spaces (Proposition \ref{productestimate}).

(iii) {\sc The Hermitian symmetric spaces $M = SO_{2k+2}/U_{k+1}$.} We first develop a general theory that applies to all irreducible Hermitian symmetric spaces. To begin with, we prove that every maximal totally geodesic submanifold $\Sigma$ with $\codim(\Sigma)$ less than half the dimension of the Hermitian symmetric space must be a complex submanifold (Proposition \ref{maxHerm}). We then prove that if the codimension of $\Sigma$ satisfy a certain inequality, then $\Sigma$ is reflective (Proposition \ref{estimate}). We can use this inequality to prove that the Index Conjecture is valid for all classical irreducible Hermitian symmetric spaces (Theorem \ref{Hermsymm}), hence in particular for the remaining space $M = SO_{2k+2}/U_{k+1}$. 

We can now state the main result of this paper: 

\begin{thm} \label{indexconj}
For every irreducible Riemannian symmetric space $M$ we have $i(M) = i_r(M)$, unless $M = G_2/SO_4$ or $M = G_2^2/SO_4$. 

Equivalently, if $M$ is an irreducible Riemannian symmetric space different from $G_2/SO_4$ and $G_2^2/SO_4$, then there exists an isometric involution $\sigma$ on $M$ so that $i(M) = \codim(\Sigma)$, where $\Sigma$ is a connected component of the fixed point set of $\sigma$ of maximal dimension. 
\end{thm}

Theorem \ref{indexconj} follows from \cite{Wo63} (for $\rk(M) = 1$), \cite{O80} (for $\rk(M) = 2$), \cite{BO16}, \cite{BO17}, \cite{BO18} (for compact simple Lie groups and many symmetric spaces of higher rank), \cite{BOR19} (for exceptional symmetric spaces), and finally Theorems \ref{LagGras}, \ref{quatGras} and \ref{Hermsymm} in this paper. In Table \ref{indexlist} we list the index for all irreducible Riemannian symmetric spaces $M$ of noncompact type together with examples of totally geodesic submanifolds $\Sigma$ with $\codim(\Sigma) = i(M)$. The symmetric spaces in the table are ordered according to their root systems ($A_r$, $B_r$, $C_r$, $D_r$, $BC_r$, $E_6$, $E_7$, $E_8$, $F_4$, $G_2$).

\begin{table}
\caption{The index $i(M)$ of irreducible Riemannian symmetric spaces $M$ of noncompact type and examples of totally geodesic submanifolds $\Sigma$ of $M$ with $\codim(\Sigma) = i(M)$} 
\label{indexlist} 
%\resizebox{12cm}{!} {
{\begin{tabular}{ | p{3.5cm}  p{4.8cm}  p{1.7cm}  p{1cm}  p{2.2cm}  |}
\hline \rule{0pt}{4mm}
\hspace{-1mm}$M$ & $\Sigma$ & $\dim(M)$ & $i(M)$ & Comments \\[1mm]
\hline \rule{0pt}{4mm}
\hspace{-2mm} 
$SO^o_{1,1+k}/SO_{1+k}$ & $SO^o_{1,k}/SO_k$ & $k+1$ & $1$ & $k \geq 1$ \\
$SL_{r+1}(\RR)/SO_{r+1}$ & $\RR \times SL_r(\RR)/SO_r$ & $\frac{1}{2}r(r+3)$ & $r$ & $r \geq 2$ \\
$SL_3(\CC)/SU_3$ & $SL_3(\RR)/SO_3$ & $8$ & $3$ & \\
$SL_4(\CC)/SU_4$ & $Sp_2(\CC)/Sp_2$ & $15$ & $5$ & \\
$SL_{r+1}(\CC)/SU_{r+1}$ & $\RR \times SL_r(\CC)/SU_r$ & $r(r+2)$ & $2r$ & $r \geq 4$ \\
$SU^*_6/Sp_3$ & $SL_3(\CC)/SU_3$ & $14$ & $6$ & \\
$SU^*_8/Sp_4$ & $Sp_{2,2}/Sp_2Sp_2$ & $27$ & $11$ & \\
$SU^*_{2r+2}/Sp_{r+1}$ & $\RR \times SU^*_{2r}/Sp_r$ & $r(2r+3)$ & $4r$ & $r \geq 4$ \\
$E_6^{-26}/F_4$ & $F_4^{-20}/Spin_9$ & $26$ & $10$ & \\[1mm]
\hline \rule{0pt}{4mm}
\hspace{-2mm} 
$SO^o_{r,r+k}/SO_{r}SO_{r+k}$ & $SO^o_{r,r+k-1}/SO_{r}SO_{r+k-1}$ & $r(r+k)$ &  $r$& $r \geq 2, k \geq 1$ \\
$SO_{2r+1}(\CC)/SO_{2r+1}$ & $SO_{2r}(\CC)/SO_{2r}$ & $r(2r+1)$ & $2r$ & $r \geq 2$ \\[1mm]
\hline \rule{0pt}{4mm}
\hspace{-2mm} 
$Sp_r(\RR)/U_r$ & $\RR H^2 \times Sp_{r-1}(\RR)/U_{r-1}$ & $r(r+1)$ & $2r-2$ & $r \geq 3$ \\
$SU_{r,r}/S(U_rU_r)$ & $SU_{r-1,r}/S(U_{r-1}U_r)$ & $2r^2$ & $2r$ & $r \geq 3$ \\
$Sp_r(\CC)/Sp_r$ & $\RR H^3 \times Sp_{r-1}(\CC)/Sp_{r-1}$ & $r(2r+1)$ & $4r-4$ & $r \geq 3$ \\
$SO^*_{4r}/U_{2r}$ & $SO^*_{4r-2}/U_{2r-1}$ & $2r(2r-1)$ & $4r-2$ &  $r \geq 3$ \\
$Sp_{2,2}/Sp_2Sp_2$ & $Sp_2(\CC)/Sp_2$ &  $16$ & $6$ &  \\
$Sp_{r,r}/Sp_rSp_r$ & $Sp_{r-1,r}/Sp_{r-1}Sp_r$ &  $4r^2$ & $4r$ & $r \geq 3$ \\
$E_7^{-25}/E_6U_1$ & $E_6^{-14}/Spin_{10}U_1$ & $54$ & $22$ & \\[1mm]
\hline \rule{0pt}{4mm}
\hspace{-2mm} $SO^o_{r,r}/SO_{r}SO_{r}$ & $SO^o_{r-1,r}/SO_{r-1}SO_{r}$ & $r^2$ & $r$ &  $r \geq 4$ \\
$SO_{2r}(\CC)/SO_{2r}$ & $SO_{2r-1}(\CC)/SO_{2r-1}$ & $r(2r-1)$ & $2r-1$ &  $r \geq 4$ \\[1mm]
\hline \rule{0pt}{4mm}
\hspace{-2mm} $SU_{r,r+k}/S(U_rU_{r+k})$ & $SU_{r,r+k-1}/S(U_rU_{r+k-1})$ & $2r(r+k)$ & $2r$ &  $r \geq 1, k \geq 1$ \\
$Sp_{r,r+k}/Sp_rSp_{r+k}$ & $Sp_{r,r+k-1}/Sp_rSp_{r+k-1}$ & $4r(r+k)$ & $4r$ &  $r \geq 1, k \geq 1$ \\
$SO^*_{4r+2}/U_{2r+1}$ &$SO^*_{4r}/U_{2r}$ & $2r(2r+1)$ & $4r$ & $r \geq 2$  \\
$F_4^{-20}/Spin_9$ & $SO^o_{1,8}/SO_8$, $Sp_{1,2}/Sp_1Sp_2$ & $16$ & $8$ &   \\
$E_6^{-14}/Spin_{10}U_1$ & $SO^*_{10}/U_5$ & $32$ & $12$ &  \\
[1mm]
\hline \rule{0pt}{4mm}
\hspace{-2mm}  $E_6^6/Sp_4$ & $F_4^4/Sp_3Sp_1$ &  $42$ & $14$ & \\
$E_6(\CC)/E_6$ & $F_4(\CC)/F_4$ & $78$ &  $26$ & \\[1mm]
\hline \rule{0pt}{4mm}
\hspace{-2mm} $E_7^7/SU_8$ & $\RR \times E^6_6/Sp_4$ & $70$ & $27$ & \\
$E_7({\mathbb C})/E_7$ & $\RR \times E_6(\CC)/E_6$ & $133$ & $54$ & \\[1mm]
\hline \rule{0pt}{4mm}
\hspace{-2mm} $E_8^8/SO_{16}$ & $\RR H^2 \times E_7^7/SU_8$ & $128$ & $56$ & \\
$E_8(\CC)/E_8$ & $\RR H^3 \times E_7(\CC)/E_7$ & $248$ & $112$ & \\[1mm]
\hline \rule{0pt}{4mm}
\hspace{-2mm} 
$F_4^4/Sp_3Sp_1$ & $SO^o_{4,5}/SO_4SO_5$ & $28$ & $8$ & \\
$E_6^2/SU_6Sp_1$ & $F_4^4/Sp_3Sp_1$ & $40$ & $12$ &  \\
$F_4(\CC)/F_4$ & $SO_9(\CC)/SO_9$ & $52$ & $16$& \\
$E_7^{-5}/SO_{12}Sp_1$ & $E_6^2/SU_6Sp_1$ & $64$ & $24$ & \\
$E_8^{-24}/E_7Sp_1$ & $E_7^{-5}/SO_{12}Sp_1$ & $112$ & $48$ & \\
[1mm]
\hline \rule{0pt}{4mm}
\hspace{-2mm} $G_2^2/SO_4$ & $SL_3(\RR)/SO_3$ & $8$ & $3$ &  \\
$G_2(\CC)/G_2$ & $G_2^2/SO_4$, $SL_3(\CC)/SU_3$ & $14$ & $6$ &\\[1mm]
\hline
\end{tabular}}
%}
\end{table}

\section{Preliminaries and notations} \label{prel}

In this section we introduce notations that we are using throughout the paper. For the general theory about Riemannian symmetric spaces we refer to \cite{H01} and \cite{Wo84}.

Let $M$ be a connected Riemannian symmetric space and $o \in M$. We normally denote by $n = \dim(M)$ the dimension of $M$ and by $r = \rk(M)$ the rank of $M$. The isometry group of $M$ is denoted by $I(M)$ and the connected component of $I(M)$ containing the identity transformation is denote by $G = I(M)^o$. We denote by $K$ the isotropy group of $G$ at $o$. Then $M$ can be identified in the canonical way with the homogeneous space $G/K$ equipped with a suitable $G$-invariant Riemannian metric. We denote $\cg$ and $\ck$ the Lie algebras of $G$ and $K$, respectively. The induced Cartan decomposition of $\cg$ is $\cg = \ck \oplus \cp$. We identify $\cp$ with the tangent space $T_oM$ of $M$ at $o$ in the usual way.

Let $\Sigma$ be a connected complete totally geodesic submanifold of $M$. Since $G$ acts transitively on $M$, we can always assume without loss of generality that $o \in \Sigma$. The tangent space $T_o\Sigma$ is a Lie triple system in $\cp$, that is, $[[T_o\Sigma,T_o\Sigma],T_o\Sigma] \subseteq T_o\Sigma$. We define $\ck' = [T_o\Sigma,T_o\Sigma] \subseteq \ck$ and $\cg' = \ck' \oplus T_o\Sigma \subseteq \cp$ and denote by $K'$ and $G'$ the connected closed subgroups of $K$ and $G$ with with Lie algebras $\ck'$ and $\cg'$, respectively. Then $\Sigma$ can be identified with the homogeneous space $G'/K'$. The group $G'$ is known as the group of glide transformations of $\Sigma$ and $K'$ as the glide isotropy group of $\Sigma$ at $o$. The normal space of $\Sigma$ at $o$ is denoted by $\nu_o\Sigma$. A Lie triple system $V$ in $\cp$ is said to be a reflective Lie triple system if the orthogonal complement of $V$ in $\cp$ is a Lie triple system. The totally geodesic submanifolds corrersponding to reflective Lie triple systems are called reflective submanifolds.

Let $M = G/K$ be a Riemannian symmetric space of compact type and consider the complexification $\cg^\CC$ of $\cg$. Using the Cartan decomposition $\cg = \ck \oplus \cp$, we define a subalgebra $\cg^*$ of $\cg^\CC$ by $\cg^* = \ck \oplus i\cp$. Let $G^*$ be the connected closed subgroup of $G^\CC$ with Lie algebra $\cg^*$. Then $G^*/K$ is a Riemannian symmetric space of noncompact type. If we start with a Riemannian symmetric space of noncompact type and perform the analogous construction, we end up with a Riemannian symmetric space of compact type. This process is known as duality between Riemannian symmetric spaces of compact type and of noncompact type. It essentially says that, up to possible finite subcoverings in the compact case, there is a one-to-one correspondence between Riemannian symmetric spaces of compact type and of noncompact type. If $V$ is a Lie triple system in $\cp$, then $iV$ is a Lie triple system in $i\cp$. Therefore, duality preserves totally geodesic submanifolds. For this reason we sometimes switch between symmetric spaces of compact type and of noncompact type, which has the advantage that we can apply methods that are specifically designed to the compact or the noncompact situation. 

\section{Reducible totally geodesic submanifolds with rank one factors} \label{tgsrk1}

In this section we investigate the codimension of locally reducible totally geodesic submanifolds with a factor of rank $1$. The following result was proved in \cite{BOR19}.

\begin{prop} \cite[Proposition 5.6]{BOR19} \label{hyperbolicfactors} 
Let $\Sigma$ be a reducible maximal totally geodesic submanifold of an irreducible Riemannian symmetric space $M$ of noncompact type. Assume that the de Rham decomposition of $\Sigma$ contains a real hyperbolic space $\RR H^k$ ($k \geq 2$), a complex hyperbolic space $\CC H^k$ ($k \geq 2$), the symmetric space $SL_3(\RR)/SO_3$, or the symmetric space $SO^o_{2,2+k}/SO_2SO_{2+k}$ ($k \geq 1$ odd). Then either $\Sigma = \RR H^{k_1} \times \RR H^{k_2}$ for some $k_1,k_2 \geq 2$, or there exists a reflective submanifold $\Sigma'$ of $M$ with $\dim(\Sigma') \geq \dim(\Sigma)$.
\end {prop}

The purpose of this section is to prove the following result:  

\begin{thm} \label{rankonefactors} 
Let $\Sigma$ be a locally reducible maximal totally geodesic submanifold of an irreducible Riemannian symmetric space $M$. Assume that the de Rham decomposition of the Riemannian universal covering space $\tilde\Sigma$ of $\Sigma$ contains a symmetric space of rank $1$ and is not equal to the Riemannian product of two spaces of nonzero constant curvature. Then $\codim(\Sigma) \geq i_r(M)$.
\end {thm}

Using duality between symmetric spaces of compact type and of noncompact type, we can restrict to the case that $M$ is of noncompact type. Then we have $\tilde\Sigma = \Sigma$. Taking into account Proposition \ref{hyperbolicfactors}, it remains to consider the quaternionic hyperbolic space $\HH H^k$ ($k \geq 2$) and the Cayley hyperbolic plane $\OO H^2$ as a possible rank $1$ factor. We start with the quaternionic case.

\begin{lm} \label{QPS} 
Let $M = \HH H^n = Sp_{1,n}/Sp_1 Sp_n$ be the $n$-dimensional quaternionic hyperbolic space and consider the reflective submanifold $\Sigma = \HH H^k = Sp_{1,k}/Sp_1 Sp_k$, $1\leq k < n$, of $M$.  Let $\rho : Sp_1Sp_k \to SO(\nu _o\Sigma)$ be the slice representation of $\Sigma$. Then the Lie algebra of $\rho(Sp_1Sp_k)$ is isomorphic to $\cs\co_3$.
\end{lm}

\begin{proof} 
The Lie algebra $\cs\cp_1 \oplus \cs\cp_k$ of the isotropy group $Sp_1Sp_k$ is linearly generated by the curvature endomorphisms $R_{u,v}$ with $u,v \in T_o\Sigma$. Hence the Lie algebra of  $\rho(Sp_1Sp_k)$ is linearly generated by the restrictions $R_{u,v}\vert_{\nu_o\Sigma} \in \cs\co(\nu_o\Sigma)$ with $u,v \in T_o\Sigma$. It is well-known that, up to a positive scalar multiple, the curvature tensor $R$ of $\HH H^n$ is given by 
\begin{equation} \label{CurQPS}
R_{u,v}w  =  - \langle v,w \rangle u + \langle u,w\rangle v 
 - \sum_{\nu=1}^3 \left(\langle J_\nu v,w \rangle J_\nu u - \langle J_\nu u,w \rangle J_\nu v- 2\langle J_\nu u,v \rangle J_\nu w \right),
\end{equation} 
where $J_1,J_2,J_3$ is a canonical basis of the quaternionic K\"{a}hler structure of $\HH H^n$ at $o$. It follows immediately from (\ref{CurQPS}) that the Lie algebra of $\rho(Sp_1Sp_k)$ is linearly generated by $J_1,J_2,J_3$ and isomorphic to $\cs\co_3$.
\end{proof}

\begin{lm} \label{QPS2} 
Let $M = \HH H^n = G/K = Sp_{1,n}/Sp_1Sp_n$, $\Sigma = \HH H^{n-1} = Sp_{1,n-1}/Sp_1Sp_{n-1}$ and $\Sigma^\perp = H/L = Sp_{1,1}/Sp_1Sp_1 \cong \HH H^1$ be the reflective submanifold of $M$ perpendicular to $\Sigma$ at $o$. Let $\tau \in I(M)$ be the geodesic reflection of $M$ in $\Sigma$. Then $\tau \in L \subset K$.
\end{lm}

\begin{proof}  
Since $\Sigma$ is a reflective submanifold of $M$, $\tau$ is an isometry. The full isometry group of $\HH H^n$ is connected and therefore $\tau \in K$. By construction, the restriction $\tau|_{\Sigma^\perp}$ is the geodesic symmetry of $\Sigma^\perp \cong \HH H^1$ at $o$. Since $\HH H^1 \cong \RR H^4$, the geodesic symmetry $\tau|_{\Sigma^\perp}$ is an inner isometry of $\Sigma^\perp$. The slice representation of $L \cong Sp_1Sp_1$ on $\nu_o\Sigma^\perp = T_o\Sigma$ is $(z,w) \cdot \xi = \xi z^{-1}$ with $(z,w) \in L \cong Sp_1Sp_1$ and $\xi \in \nu_o\Sigma^\perp$. The isotropy representation of $L \cong Sp_1Sp_1$ on $T_o\Sigma^\perp = \nu_o\Sigma$ is $(z,w) \cdot X = wXz^{-1}$ with $(z,w) \in L \cong Sp_1Sp_1$ and $X \in T_o\Sigma^\perp$. As $\tau|_\Sigma = \id_\Sigma$, we have $d_o\tau|_{T_o\Sigma} = \id|_{T_o\Sigma}$ and $d_o\tau|_{T_o\Sigma^\perp} = -\id|_{T_o\Sigma^\perp}$, it follows that $\tau \in L$, corresponding to the element $(1,-1) \in Sp_1Sp_1$.
\end{proof}

We now prove Theorem \ref{rankonefactors} for the case that $\Sigma$ contains a quaternionic hyperbolic space as a de Rham factor.

\begin {prop} \label{quat}
Let $\Sigma$ be a reducible maximal totally geodesic submanifold of an irreducible Riemannian symmetric space $M$ of noncompact type. Assume that the de Rham decomposition of $\Sigma$ contains a quaternionic hyperbolic space $\HH H^k$, $k \geq 2$. Then $\codim(\Sigma) \geq i_r(M)$.
\end {prop}

\begin {proof}
We can assume that $o \in \Sigma$ and write $M = G/K$ and $\Sigma = G'/K'$ as in Section \ref{prel}. If $\Sigma$ is nonsemisimple, then $\Sigma$ is reflective by \cite[Theorem 1.2]{BO16} and hence $\codim(\Sigma) \geq i_r(M)$. We therefore can assume that $\Sigma$ is semisimple. By assumption, we have $\Sigma = \Sigma_1 \times \Sigma_2$ with $\Sigma_1 = G'_1/K'_1 \cong \HH H^k = Sp_{1,k}/Sp_1Sp_k$. We fix a totally geodesic $\tilde\Sigma \cong \HH H^{k-1}$ in $\Sigma_1$ with $o \in \tilde\Sigma$. This is a reflective submanifold of $\Sigma_1$ and there exists a reflective submanifold $\tilde\Sigma^\perp \cong \HH H^1$ of $\Sigma_1$ with $o \in \tilde\Sigma^\perp$ that is perpendicular to $\tilde\Sigma$ at $o$. Let $\tau \in I(\Sigma_1)$ be the geodesic reflection of $\Sigma_1$ in $\tilde\Sigma$. By Lemma \ref{QPS2}, $\tau \in K'_1 \subset K' \subset K$. With the same arguments as in the proof of Proposition \ref{hyperbolicfactors}, by replacing $\tau$ with a suitable odd power of $\tau$, we may assume that $\tau$ is an involutive isometry of $M$. Moreover, analogously to the proof of Proposition \ref{hyperbolicfactors}, if the set 
\[
V = \Fix_{\nu_o\Sigma}(d_o\tau) =  \{v \in \nu_o\Sigma : d_o\tau(v) = v\}
\] 
of fixed vectors of $d_o\tau$ in $\nu_o\Sigma$ is trivial, then $\tilde\Sigma \times \Sigma_2 \subset \Sigma$ is a reflective submanifold of $M$, and it follows from \cite[Corollary 2.9]{BOR19} that $\Sigma$ is a reflective submanifold of $M$, which implies $\codim(\Sigma) \geq i_r(M)$.

Thus we can assume $\dim(V) \geq 1$. Since $\tau$ is involutive, the totally geodesic submanifold $\Sigma'$ of $M$ with $T_o\Sigma' = T_o\tilde\Sigma \oplus T_o\Sigma_2 \oplus  V$ is reflective. If $\dim(V) \geq 4$, then $\dim(\Sigma) \leq \dim(\Sigma')$ and thus $\codim(\Sigma) \geq i_r(M)$. If $\dim(V) =1$, we obtain by a similar argument to that used in the proof of Proposition \ref{hyperbolicfactors}, that $\Sigma$ has only one other de Rham factor, which is isometric to a real hyperbolic space. Then, again by Proposition \ref{hyperbolicfactors}, $\Sigma$ is a product of real hyperbolic spaces, which is a contradiction. 
Thus we are left with the two possibilities $\dim(V) \in \{2, 3\}$. 

\begin{lm} \label{lemmaA}
$T_o\Sigma_2 \oplus V$ is a Lie triple system.
\end{lm}
 
\begin{proof}
Consider the slice representation $\bar\rho$ of the isotropy group $Sp_1Sp_{k-1}$ of $\tilde\Sigma$ on its normal space $T_o\Sigma_2 \oplus V$ in $\Sigma'$. The isotropy group $Sp_1Sp_{k-1}$ acts trivially on $T_o\Sigma_2$, since $\Sigma_1 \times \Sigma_2 = \Sigma$ and $\tilde\Sigma \subset \Sigma_1$. Since $\dim(Sp_1Sp_{k-1}) \geq 6$ and $\dim(SO(V)) \leq 3$, $\ker(\bar\rho)$ is a nontrivial normal subgroup of $Sp_1Sp_{k-1}$. The set $F$ of fixed vectors of this normal subgroup in $T_o\tilde\Sigma$ is trivial. (Note that $F \neq T_o\tilde\Sigma$, since $Sp_1Sp_{k-1}$ acts  almost effectively.) In fact, if $F$ is nontrivial, $F$ must be invariant under $Sp_1Sp_{k-1}$ and so $\tilde\Sigma$ would be reducible, which is a contradiction. Therefore the set of fixed vectors of $\ker(\bar \rho)$ in 
$T_o\Sigma'$ is exactly $T_o\Sigma_2 \oplus V$ and so this subspace is a Lie triple system. 
\end{proof}

\begin{lm} \label{lemmaB}
The slice representation $\bar\rho$ of the isotropy group $Sp_1Sp_{k-1}$ of $\tilde\Sigma$ on its normal space $T_o\Sigma_2 \oplus V$ in $\Sigma'$ is trivial.
\end{lm}

\begin{proof}
Since $Sp_1Sp_{k-1}$ acts trivially on $T_o\Sigma_2$, we only need to show that $Sp_1Sp_{k-1}$ acts trivially on $V$. Assume that $\bar\rho$ is nontrivial. Since $Sp_1Sp_{k-1}$ has no normal subgroups of codimension $1$, we must have $\dim(\bar\rho(Sp_1Sp_{k-1})) > 1$. This proves our assertion for $\dim(V) = 2$. Assume that $\dim(V)=3$ and $\bar\rho(Sp_1Sp_{k-1}) = SO(V)$. Then there are no nonzero fixed vectors by $\bar\rho(Sp_1Sp_{k-1})$ in $V$. From Lemma \ref{lemmaA} we know that  $\tilde\Sigma$ is a reflective submanifold of $\Sigma'$. So the set of fixed vectors of $\bar\rho$, which coincides with $T_o\Sigma_2$, is invariant under the isotropy group of the perpendicular reflective submanifold $P$, where $T_oP = T_o\Sigma_2 \oplus V$. Then $\Sigma_2$ is a de Rham factor of $P$ and so $[T_o\Sigma_2 , V]=\{0\}$. This implies that the centralizer $\cz_\cp(T_o\Sigma_2)$ of $T_o\Sigma_2$ in $\cp$ contains $T_o\Sigma_1 \oplus V$. Then $\cz_\cp(T_o\Sigma_2) + T_o\Sigma_2$ is a proper Lie triple system in $T_oM$ containing  $T_o\Sigma$ properly. This contradicts the maximality of $\Sigma$ and so the assertion follows.
\end{proof}

\begin{lm} \label{lemmaC}
$\tilde\Sigma$ is a de Rham factor of $\Sigma'$.
\end{lm}

\begin{proof}
Since $\bar\rho$ is trivial by Lemma \ref{lemmaB}, it follows from \cite[Proposition 3.8]{BOR19} that $\tilde\Sigma$ is either a de Rham factor of $\Sigma'$ or it is contained in a de Rham factor $M_1$ of $\Sigma'$ of constant curvature. Assume the latter holds. Note that $M_1$ is strictly contained in $\Sigma'$. In fact, we have $\rk(\Sigma') \geq 2$ since  $[T_o\tilde\Sigma, T_o\Sigma_2] = \{0\}$. Let us write $\Sigma' = M_1\times M_2$, where $M_2$ is not necessarily irreducible. Note that $\Sigma_2 \subset M_2$. In fact, if $v \in T_o\Sigma_2$, then its orthogonal projection onto $T_oM_1$ must be trivial since $\rk(M_1) = 1$, $\tilde\Sigma \subset M_1$ and $[v, T_o\tilde\Sigma] = \{0\}$. Then $[T_o\Sigma_2, T_oM_1] = \{0\}$ and, as in the proof of Lemma \ref{lemmaB}, $\cz_\cp(T_o\Sigma_2) + T_o\Sigma_2$ is a proper Lie triple system in $T_oM$ that properly contains $T_o\Sigma$. This contradicts the maximality of $\Sigma$ and the assertion follows.
\end{proof}

We continue with the proof of Proposition \ref{quat}. Let $\Sigma'^\perp$ be the reflective submanifold of $M$ which is perpendicular to $\Sigma'$ at $o$. Then, by construction, $T_o\Sigma'^\perp$ is the $(-1)$-eigenspace of $d_o\tau$. Note that $T_o\tilde\Sigma^\perp \subset T_o\Sigma'^\perp$. According to Lemma \ref{QPS}, the kernel $H$ of the representation of $Sp_1Sp_{k-1}$ (as in Lemma \ref{lemmaB}) on $T_o\tilde\Sigma^\perp$ must be isomorphic to  $Sp_{k-1}$ (and hence of dimension $\geq 3$). It follows from Lemma \ref{lemmaC} that $H$ acts trivially on $T_o\Sigma_2 \oplus V$. Since $H$ is a normal subgroup of the isotropy group at $o$ of the de Rham factor $\tilde\Sigma$ of $\Sigma'$ (see Lemma \ref{lemmaC}), $H$ is a normal subgroup of the isotropy group of $\Sigma'$. Then, the set $W$ of fixed vectors of $H$ on the normal space $\nu _o\Sigma'$, which contains $T_o\tilde\Sigma^\perp$, is invariant under the isotropy group at $o$ of the complementary reflective submanifold $\Sigma'^\perp$. Then $W =T_oQ$, where $Q$ is a de Rham factor of $\Sigma'^\perp$. 

Assume that $W = T_o\Sigma'^\perp$. Using Lemma \ref{lemmaB} we see that $H$ acts trivially on $T_o\Sigma_2 \oplus V$. Then $H$ acts trivially on $\nu _o\tilde\Sigma = T_o\Sigma_2 \oplus V \oplus W$. Using \cite[Proposition 2.8]{BOR19} we obtain that $\tilde\Sigma$ is reflective. Then, using \cite[Corollary 2.9]{BOR19} and the fact that  $\tilde\Sigma \subset \Sigma$, we see that $\Sigma$ is reflective and hence $\codim(\Sigma) \geq i_r(M)$. 

Next, assume that $W$ is a proper subspace of $T_o\Sigma'^\perp$, or equivalently, $Q$ is properly contained in $\Sigma'^\perp$. Let us write, as a nontrivial Riemannian product, $\Sigma'^\perp = Q \times Q'$, where  $o \in Q'$ and $T_oQ'$ is the orthogonal complement of $W$ in $T_o\Sigma'^\perp$. Then $\dim(Q') \geq 3$, because otherwise $H \cong Sp_{k-1}$ would act trivially on $T_oQ'$. Let us consider the involutive isometry $\tau \in L$, where $L \cong Sp_1Sp_1\subset K$ is the glide isotropy group of $\tilde\Sigma^\perp$. Recall that $T_o\Sigma'$ is the $(+1)$-eigenspace of $d_o\tau$ and $T_o\tilde\Sigma^\perp$ is contained in the $(-1)$-eigenspace of $d_o\tau$. 

\begin{lm} \label{lemmaD}
The involution $\tau$ commutes with every isometry in the glide isotropy group $L \cong Sp_1Sp_1$ of $\tilde\Sigma^\perp$.
\end{lm}

\begin{proof}
Consider the set  $A =\{k\circ \tau \circ k^{-1} \circ \tau^{-1}: k\in L\} \subset L$ of isometries of $M$. Note that $A$ is connected, since $L$ is connected. Note also that any isometry in $A$ acts trivially on $\tilde\Sigma^\perp$. Since $L$ acts almost effectively on $\tilde\Sigma^\perp$, $A$ must be discrete and thus $A = \{\id_M\}$, since $A$ is connected. 
\end{proof}

From Lemma \ref{lemmaD} we see that $L$, via the isotropy representation, leaves the $(+1)$-eigenspace $T_o\Sigma'$ of $d_o\tau$ invariant. Consequently, $L$ leaves $\Sigma'$ invariant. According to Lemma \ref{QPS} there exists a nontrivial normal subgroup $\tilde L \cong Sp_1$ of $L$ such that $\tilde L$ acts trivially on $T_o\tilde\Sigma$. Note that $\tilde L$ acts trivially also on $T_o\Sigma_2$, which follows from the fact that $\Sigma = \Sigma_1 \times \Sigma_2$ (note that $\tilde\Sigma^\perp$ is a totally geodesic submanifold of $\Sigma_1$).

\begin{lm} \label{lemmaE}
The normal subgroup $\tilde L$ of $L$ acts trivially on $V$.
\end{lm}

\begin{proof}
We have $\tilde L (V) \subset V$, because $\tilde L (T_o\Sigma') = T_o\Sigma' = T_o\tilde\Sigma \oplus T_o\Sigma_2 \oplus V$ and $\tilde L$ acts trivially on both $T_o\tilde\Sigma$ and $T_o\Sigma_2$. If $\dim(V)= 2$, then $\tilde L$ acts trivially on $V$ since $\tilde L \cong Sp_1 \cong Spin_3$. Let $\dim(V)= 3$. Assume that $\tilde L$ acts on $V$ nontrivially and let $g \in \tilde L$ be such that $h = g|_V \neq \id$. Since $\tilde L$ is connected, $+1$ is an eigenvalue of $d_oh$ with multiplicity $1$. Let $\RR v$ with $0 \neq v \in V$ be the corresponding eigenspace. Recall from Lemma \ref{lemmaA} that  $T_o\Sigma_2 \oplus V$ is a Lie triple system. This Lie triple system is invariant under $\tilde L$. Then $T_o\Sigma_2 \oplus \RR v$ is also a Lie triple system, since it coincides with the set of fixed vectors of $g$ in $T_o \Sigma_2 \oplus V$. Let $X$ be the totally geodesic submanifold of $M$ with $T_oX = T_o\Sigma_2 \oplus \RR v$. Then $\Sigma_2$ is a semisimple  totally geodesic hypersurface of $X$. Then, by \cite[Lemma 5.5]{BOR19}, either there exists an irreducible de Rham factor $\Sigma_2'$ of $\Sigma_2$ with constant curvature, or  $X$ is a Riemannian product $X = \Sigma_2 \times \RR$. In the first case, by Proposition \ref{hyperbolicfactors}, $\Sigma$ is a product of spaces of constant curvature, which contradicts our assumption. In the second case, $\cz_\cp(T_o\Sigma_2') + T_o\Sigma_2'$ is a proper Lie triple system in $T_oM$ that properly contains $T_o\Sigma$, since it also contains $v$. This contradicts the maximality of $\Sigma$.  
\end{proof}
 
From Lemma \ref{lemmaE} and its preceding paragraph we see that $\tilde L$ acts trivially on $T_o\Sigma'$. Since $\tilde L \subset L$, and $L$ is included in the glide isotropy group of $Q$ at $o$, $\tilde L$ acts trivially on $T_oQ'$. Let $U \subset T_oM$ be the subspace of fixed vectors of $\tilde L$. Then $T_o\Sigma' \oplus T_oQ' \subset U$ and so
 \[
 \dim (U) \geq \dim (T_o\Sigma') + \dim (T_oQ') \geq \dim (\Sigma) -2 + 3 > \dim(\Sigma).
 \]
Then the totally geodesic submanifold $S$ of $M$ with $T_oS =U$ satisfies $\dim(U) > \dim(\Sigma)$.  Moreover, $S$ contains the reflective submanifold $\Sigma'$. Then $S$ is reflective by \cite[Corollary 2.9]{BOR19} and it follows that $\codim(\Sigma) > \codim(S) \geq i_r(M)$. This finishes the proof of Proposition \ref{quat}. 
\end{proof}

We now consider the Cayley hyperbolic plane $\OO H^2$ as a possible rank $1$ factor. 

\begin{prop} \label{Cayley}
Let $\Sigma$ be a reducible maximal totally geodesic submanifold of an irreducible Riemannian symmetric space $M = G/K$ of noncompact type. Assume that the de Rham decomposition of $\Sigma$ contains the Cayley hyperbolic plane $\OO H^2$ as a factor. Then $\codim(\Sigma) \geq i_r(M)$.
\end{prop}

\begin{proof}
The full isometry group of $\OO H^2$ is connected and isomorphic to the noncompact real simple Lie group $F_4^{-20}$. The isotropy group at $o$ is isomorphic to $Spin_9$ and thus we can write $\OO H^2 = F_4^{-20}/Spin_9$. The isotropy representation of $Spin_9$ on $T_o\OO H^2$ is equivalent to the spin representation of $Spin_9$ on $\RR^{16}$.

As usual, we can assume $o \in \Sigma$ and write $\Sigma = G'/K'$ as in Section \ref{prel}. If $\Sigma$ is nonsemisimple, then $\Sigma$ is reflective by \cite[Theorem 1.2]{BO16} and thus $\codim(\Sigma) \geq i_r(M)$. We can therefore assume that $\Sigma$ is semisimple. By assumption, we have $\Sigma = \OO H^2 \times \bar\Sigma$, where $\bar\Sigma$ is a semisimple Riemannian symmetric space of noncompact type. 

The Cayley hyperbolic plane admits only one type of polars, namely Cayley hyperbolic lines $\OO H^1$, all of which are congruent to each other in $\OO H^2$ and isometric to the real hyperbolic space $\RR H^8$. We choose a Cayley hyperbolic line $P \cong \OO H^1$ in $\OO H^2$ with $o \in P$ and denote by $Q \cong \OO H^1$ the Cayley hyperbolic line in $\OO H^2$ with $o \in Q$ that is perpendicular to $P$ at $o$. Note that $P$ and $Q$ is a pair of complementary reflective submanifolds of $\OO H^2$ and $T_o\OO H^2 = T_oP \oplus T_oQ$. The subgroup of the isotropy group $Spin_9$ leaving this decomposition invariant is (isomorphic to) $Spin_8$. The restriction to $Spin_8$ of the isotropy representation of $Spin_9$ is equivalent to the direct sum of the two inequivalent spin representations of $Spin_8$ on $\RR^8$. The subgroup $Spin_8$ is the isotropy group of each of the two groups of glide transformations of $P$ and $Q$.

We denote by $\tau_P,\tau_Q \in Spin_8$ the geodesic reflections of $\OO H^2$ in $P$ and $Q$, respectively. Since $Spin_8 \subset Spin_9 \subset K' \subset K$, both $\tau_P$ and $\tau_Q$ can be viewed as isometries of $M$. Note that $\tau_P$ and $\tau_Q$ lie both in the (finite) center $Z_{Spin(8)}$ of $Spin_8$, since $\tau_P|_Q, \tau _Q|_P$ are the geodesic symmetries of $Q$ and $P$, respectively (and $Spin_8$ acts almost effectively on both $P$ and $Q$). Since any nontrivial element in the center $Z_{Spin(8)}$ of $Spin _8$ has order $2$, $\tau _P$ and $\tau _Q$ have both order $2$, as elements of $K$. We define $V_P = \Fix_{\nu_o\Sigma}(d_o\tau_P) = \{ v \in \nu_o\Sigma : d_o\tau_P(v) = v\}$. Then $T_oP \oplus T_o\bar{\Sigma} \oplus V_P = \Fix_{T_oM}(d_o\tau_P)$. Evidently, $\Fix_{T_oM}(d_o\tau_P)$ is a reflective Lie triple system in $T_oM$. Let $\Sigma_P$ be the reflective submanifold of $M$ with $T_o\Sigma_P = \Fix_{T_oM}(d_o\tau_P)$ and $\Sigma_P^\perp$ be the reflective submanifold of $M$ with $T_o\Sigma_P^\perp = \nu_o\Sigma_P$. Note that $Q \subseteq \Sigma_P^\perp$.

If $\dim(V_P) = 0$, then $P \times \bar{\Sigma}$ is a reflective submanifold of $M$. It then follows from \cite[Corollary 2.9]{BOR19} that $\Sigma$ is a reflective submanifold of $M$ and hence $\codim(\Sigma) \geq i_r(M)$. 

If $\dim(V_P) \geq 8$, then $\dim(\Sigma) \leq \dim(\Sigma_P)$ and hence $\codim(\Sigma) \geq \codim(\Sigma_P) \geq i_r(M)$. 

If $\dim(V_P) \in \{1,\ldots,7\}$, then the isotropy group $Spin_8$ acts trivially on $V_P$. Since $Spin_8$ acts trivially also on $T_o\bar\Sigma$, it follows that the slice representation of $Spin_8$ on the normal space $T_o\bar{\Sigma} \oplus V_P$ of $P$ at $o$ in $\Sigma_P$ is trivial. An analogous argument as for the quaternionic case in the proof of Lemma \ref{lemmaC} shows that $P$ is a de Rham factor of $\Sigma_P$. The set $\Fix_{T_o\Sigma_P^\perp}(d_o\tau_Q)$ is a Lie triple system containing $T_oQ$ and invariant under the glide isotropy group $H$ of $\Sigma_P^\perp$ at $o$. In fact, $H$ leaves invariant the factor $P$ of $\Sigma_P$ and so it must leave invariant the finite center $Z_{Spin(8)}$ of $Spin_8$. Then the identity component $H^o$ of $H$ must commute with $Z_{Spin(8)}$ and in particular with $\tau_Q$. Thus there exists a Riemannian factor $\tilde\Sigma_P^\perp$ of $\Sigma_P^\perp$ such that $T_o\tilde\Sigma_P^\perp = \Fix_{T_o\Sigma_P^\perp}(d_o\tau_Q)$. Note that $Q \subseteq \tilde\Sigma_P^\perp$. We have $d_o\tau_Q|_{T_oP} = -\id|_{T_oP}$ and, since $\tau_Q \in Spin_8$ and the slice representation of $Spin_8$ on the normal space $T_o\bar\Sigma \oplus V_P$ of $P$ at $o$ in $\Sigma_P$ is trivial, $d_o\tau_Q|_{T_o\bar\Sigma \oplus V_P} = \id|_{T_o\bar\Sigma \oplus V_P}$. Therefore, if $\tilde\Sigma_P^\perp = \Sigma_P^\perp$, then $P$ is a reflective submanifold of $M$ and  \cite[Corollary 2.9]{BOR19} implies that $\Sigma$ is a reflective submanifold of $M$ and hence $\codim(\Sigma) \geq i_r(M)$. Otherwise, we get a proper Riemannian product decomposition $\Sigma_P^\perp = \tilde\Sigma_P^\perp \times \bar\Sigma_P^\perp$. Since $Q \subset \tilde\Sigma_P^\perp$, the isotropy group $Spin_8$ of $Q$ acts trivially on $ T_o\bar\Sigma_P^\perp$. As $\tau_Q \in Spin_8$, this implies $d_o\tau_Q|_{T_o\bar\Sigma_P^\perp} = \id_{T_o\bar\Sigma_P^\perp}$, which is a contradiction to $T_o\tilde\Sigma_P^\perp = \Fix_{T_o\Sigma_P^\perp}(d_o\tau_Q)$. This finishes the proof.
\end{proof}

\section{General structure results} \label{gsr}

In this section we prove some general results about totally geodesic submanifolds in symmetric spaces, which will be useful for later purposes. We start by investigating reflection hyperplanes of totally geodesic submanifolds.

\begin{prop} \label{cardinal}
Let $M=G/K$ be a simply connected irreducible Riemannian symmetric space and $\Sigma = G'/K'$ be a totally geodesic submanifold. Let $W$ and $W'$ be the Weyl groups associated with $M$ and $\Sigma$, respectively. Let $b$ and $b'$ be the number of reflection hyperplanes of $W$ and $W'$, respectively. Then $b' \leq b$. 
\end{prop}

\begin{proof} 
By duality, we can assume that $M$ is of noncompact type. Then $\Sigma$ is simply connected and hence $K'$ is connected. Let $\cg = \ck \oplus \cp$ and $\cg' = \ck' \oplus \cp'$ be the Cartan decompositions associated with $(G,K)$ and $(G',K')$, respectively. As usual, we identify $T_oM$ with $\cp$ and $T_o\Sigma$ with $\cp'$. Let $\ca'$ be a maximal abelian subspace of $\cp'$ and $\ca$ be a maximal abelian subspace of $\cp$ with $\ca' \subseteq \ca$. We consider $W$ and $W'$ as reflection groups of $\ca$ and $\ca'$, respectively. 

Let $H'_1,\ldots,H'_{b'} \subset \ca'$ and $H_1,\ldots ,H_b \subset \ca$ be the distinct reflection hyperplanes associated with $W'$ and $W$, respectively. We define the set
\[
J = \{j \in \{1,\ldots,b\} :  \ca' \subseteq H_j\},
\]
which could be an empty or a nonempty set. Since the intersection of all reflection hyperplanes is $\{0\}$, $J$ is properly contained in $\{1, \ldots , b\}$. By a suitable labelling of the reflection hyperplanes we can assume that there exists $j_0\in \{1,\ldots,b\}$ such that 
\[
j \notin J \iff j\geq j_0. 
\] 
Then $H_j \cap \ca' = \ca'$ if $j < j_0$ and $H_j \cap \ca'$ is a hyperplane of $\ca'$ if $j \geq j_0$. Note that any two such hyperplanes $H_j \cap \ca'$ may coincide. 

Let us assume that there exists a reflection hyperplane $H'_d \subset \ca'$ such that for every reflection hyperplane $H_j \subset \ca$ we have $H_d' \neq H_j \cap \ca'$. This is always true for $j < j_0$, and for $j \geq j_0$ this means that $H_d' \cap H_j$ is a hyperplane of $H'_d$. Note that $H_d' \cap H_i'$ is also a hyperplane of $H'_d$ for all $i \neq d$. Therefore we can find $0 \neq u \in H'_d$ so that $u \notin H_i'$ for all $i \neq d$ and $u \notin H_j$ for all $j \geq j_0$. There exists $\epsilon > 0$ so that the open ball $B_\epsilon(u)$ in $\ca'$ with radius $\epsilon$ and center $u$ does not intersect $H_i'$ for all $i \neq d$ and $H_j$ for all $j \geq j_0$. 

We now choose a point $v \in B_\epsilon(u) \subset \ca'$ that is not contained in $H_d'$ and define the curve $\gamma : [0,1] \to B_\epsilon(u),\ t \mapsto v + t(u-v)$, which parametrizes the line segment from $v$ to $u$. By construction, we have $\{ j \in J : \gamma(t) \in H_j\} = J$
for each $t \in [0,1]$. According to the Slice Theorem of Hsiang, Palais and Terng (see \cite[Section 2]{HPT88} and \cite[Section 6.5]{PT88}), the dimension of the isotropy orbit $K \cdot \gamma(t)$ satisfies
\[
\dim (K\cdot \gamma (t)) = m - \sum_{j\in J} m_j, 
\]
where $m$ is the dimension of a principal $K$-orbit in $\cp$ and $m_j$ is the multiplicity of a focal point in $H_j$ which is not in any other reflection hyperplane $H_\nu$ for $\nu \neq j$. It follows that $\dim(K \cdot \gamma(t))$ is independent of the choice of $t \in [0,1]$.

On the one hand, this implies that the identity components of the isotropy groups $K_{\gamma (t)}$ ($t \in [0,1]$) coincide, or equivalently, the isotropy algebras $\ck_{\gamma (t)}$ ($t \in [0,1]$) coincide. On the other hand, by the choice of $u$ and $v$, the orbit $K' \cdot v \subset \cp'$ is a principal orbit of the isotropy action of $K'$ on $\cp'$ and $K' \cdot u$ is a parallel focal orbit of $K' \cdot v$. Therefore the isotropy algebra $\ck'_v$ is strictly contained in the isotropy algebra $\ck'_u$. Thus there exists $z \in \ck'$ such that $\ad(z)u = 0$ and $\ad(z)v \neq 0$. Since $\ck' \subset \ck$, this implies $z \in \ck_u$ and $z \notin \ck _v$, which contradicts $\ck_u = \ck_v$. It follows that for every reflection hyperplane $H'_i \subset \ca'$ there exists a reflection hyperplane $H_j \subset \ca$ such that $H_i' = H_j \cap \ca'$, which implies $b' \leq b$. 
\end{proof}

Let $M = G/K$ be a simply connected irreducible Riemannian symmetric space and consider the marked Dynkin diagram associated with $G/K$, which is the Dynkin diagram associated with $G/K$ together with the dimensions of the corresponding root spaces. The dimension of a root space is also called the multiplicity of the root. Here we make the convention that the multiplicity of a non-reduced root $\alpha$ is obtained by adding up the dimensions of the root spaces of $\alpha$ and $2\alpha$. We denote by $\Phi$ the corresponding root system and by $\Phi^+$ the positive roots.

If $\Phi$ is reduced and all roots have the same length, we call all roots long. If $\Phi$ is reduced and there are roots of different length, then there are exactly two different lengths and we can naturally distinguish between long and short roots. If $\Phi$ is non-reduced, we call the non-reduced roots short and the other roots long. 

Recall that the Weyl group acts transitively on the sets of long and short roots and so the multiplicities of any two long (resp.\ short) roots are the same. This implies that all long roots have the same multiplicity $m_1$, and all short roots have the same multiplicity $m_2$. If there are no short roots, our convention is $m_1 = m_2$.

We denote by $\bar{l}$ (resp.\ $\bar{s}$) the number of positive long (resp.\ short) roots in $\Phi^+$. Then we have 
\begin{equation}\label{dim3}
\dim(M) =  m_1 \bar l +  m_2 \bar s + \rk(M).
\end{equation}
We call $m_1$ and $m_2$ the {\it associated multiplicities} of $M$.

\begin{prop}\label{multcomp}
Let $M = G/K$ be an irreducible Riemannian symmetric space of noncompact type with associated multiplicities $m_1$ and $m_2$. Let $\Sigma = G'/K'$ be a totally geodesic submanifold of $M$ with $\rk(\Sigma) = \rk(M)$. Let $\Sigma_1$ be a de Rham factor of $\Sigma$ with associated multiplicites $m'_1$ and $m'_2$. Then we have
\[
\max\{m'_1, m'_2\}\leq \max\{m_1, m_2\}.
\]
\end{prop}

\begin{proof}
Let $\ca$ be a maximal abelian subspace of $\cp' \cong T_o\Sigma$. Since $\rk(\Sigma) = \rk(M)$, $\ca$ is also a maximal abelian subspace of $\cp \cong T_oM$. We choose $v \in \ca$ so that $K' \cdot v$ is a principal orbit of the $K'$-action on $\cp'$ and $K \cdot v$ is a principal orbit of the $K$-action on $\cp$. Every normal vector $\xi \in \nu_v(K\cdot v)$ of $K \cdot v$ at $v$ extends uniquely to a $K$-invariant normal vector field $\tilde\xi$ of $K\cdot v$. The restriction $\bar \xi = \tilde \xi|_{K'\cdot v}$ of $\tilde\xi$ to $K' \cdot v$ is a $K'$-invariant normal vector field of $K'\cdot v$. The actions of $K$ on $\cp$ and of $K'$ on $\cp'$ are polar (see \cite{BCO16,PT88}). This implies that $\tilde\xi$ is parallel with respect to the normal connection of $K\cdot v$ in $\cp$ and $\bar \xi$ is parallel with respect to the normal connection of $K'\cdot v$ in $\cp'$. 

Let $A$ (resp.\ $A'$) be the shape operator of $K\cdot v$ in $\cp$ (resp.\ of  $K'\cdot v$ in $\cp'$). For $X \in \ck'$ we have
\begin{equation}\label{shape2}
A_\xi (\ad(X)v) = -\left.\textstyle{\frac{d}{dt}}\right|_{t=0} \tilde\xi_{c_v(t)} 
= -\left.\textstyle{\frac{d}{dt}}\right|_{t=0} \bar\xi_{c_v(t)}
=  A'_\xi (\ad(X)v), 
\end{equation}
where $c_v(t) = \Exp(tX)v \in K' \cdot v \subset K \cdot v$. This shows that the tangent space $T_v(K'\cdot v)$ of $K' \cdot v$ at $v$ is invariant under the shape operator $A_\xi$ of $K\cdot v$ with respect to $\xi$. Therefore, each common eigenspace of the (commuting) family of shape operators of $K'\cdot v$ at $v$ is contained in a common eigenspace of the (commuting) family of shape operators of $K\cdot v$ at $v$. Moreover, any curvature normal of $K' \cdot v$ at $v$ is a curvature normal of $K \cdot v$ at $v$. According to \cite [page 63]{BCO16}, the common eigenspaces of the family of shape operators at $v$ of the principal orbit $K \cdot v$ in $\cp$ are given by $E_{\alpha} = \cp_\alpha$ for $\alpha \in \Phi^+$ reduced and $E_{\alpha} = \cp_\alpha \oplus \cp_{2\alpha}$ for $\alpha \in \Phi^+$ with $2\alpha \in \Phi^+$, where $\cp_\alpha = (\cg_\alpha \oplus \cg_{-\alpha}) \cap \cp$. Then
\begin{equation}\label {dimE3}
\dim (E_\alpha) = m_i
\end{equation}
with $i=1$ or $i=2$, depending on whether $\alpha$ is a long root or a short root. The analogous statement holds for the common eigenspaces of the family of shape operators at $v$ of the principal orbit $K' \cdot v$ in $\cp'$. This finishes the proof.
\end{proof}

Recall that a Riemannian symmetric space $M=G/K$ of compact type is inner if the geodesic symmetry $\sigma_o$ of $M$ at $o$ belongs to $K$. It is known (see e.g.\ \cite[Theorem 8.6.7]{Wo84}) that $G/K$ is inner if and only if $\rk(G) = \rk(K)$. A non-inner symmetric space is also called an outer symmetric space.

\begin {prop} \label {innSigma}
Let $M=G/K$ be an outer irreducible Riemannian symmetric space of compact type and let $\Sigma = G'/K'$ be a maximal totally geodesic submanifold of $M$. If $\Sigma$ is an inner symmetric space, then $\Sigma$ is a reflective submanifold.
\end{prop}

\begin{proof}
Let $\tau \in K'$ be the geodesic symmetry of $\Sigma$ at $o$. Using the same arguments as in the proof of Proposition \ref{hyperbolicfactors}, we may assume that $\tau$ is an involutive isometry of $M$. Since $M$ is an outer symmetric space, $\tau$ must be different from the geodesic symmetry $\sigma _o$ of $M$. Let $\tilde\Sigma$ be the connected component containing $o$ of the fixed point set of the involutive isometry $\sigma_o \circ \tau$. Then $\tilde\Sigma$ is a reflective totally geodesic submanifold of $M$ containing $\Sigma$. 
Since $\Sigma$ is maximal, $\tilde\Sigma = \Sigma$ and thus $\Sigma$ is reflective. 
\end{proof}

For details on the following constructions we refer to \cite{NT00} and the references therein. Let $M = G/K$ be a Riemannian symmetric space of compact type. The point $o$ is an isolated fixed point of the geodesic symmetry $\sigma_o$ of $M$ at $o$. The connected components different from $\{o\}$ of the fixed point set $\Fix_M(\sigma_o) = \{p \in M : \sigma_o(p) = p\}$  are so-called polars of $M$. A polar consisting of a single point is also called a pole of $o$. Every antipodal point on a closed geodesic through $o$ lies on a polar. More precisely, if $\gamma : [0,1] \to M$ is a geodesic with $\gamma(0) = o = \gamma(1)$, then $\gamma(\frac{1}{2})$ lies in a polar of $o$. In fact, the set $\Fix_M(\sigma_o) \setminus \{o\}$ coincides with the set of antipodal points of $o$. Let $o \neq p \in \Fix_M(\sigma_o)$ and denote by $M^+(p)$ the polar containing $p$. Then $M^+(p) = K \cdot p$ (even is $K$ is not connected, the orbit $K \cdot p$ is connected). Thus the polars of $o$ are the orbits of the isotropy group $K$ at $o$ through the antipodal points of $o$ 

Every polar $M^+(p)$ is a reflective submanifold of $M$. In fact, we have $T_pM^+(p) = \{v \in T_pM : d_p\sigma_o(v) = v\}$ and $\nu_pM^+(p) = \{v \in T_pM : d_p\sigma_o(v) = -v\}$. It follows that $\nu_pM^+(p)$ is the fixed point set of $d_p(\sigma_o \circ \sigma_p) = d_p(\sigma_p \circ \sigma_o)$ and therefore the connected component of $\Fix_M(\sigma_o \circ \sigma_p)$ containing $p$ is a totally geodesic submanifold of $M$ and $T_pM = T_pM^+(p) \oplus T_pM^-(p)$. It follows that $M^+(p)$ and $M^-(p)$ is a complementary pair of reflective submanifolds. Any such submanifold $M^-(p)$ is called a meridian. Any meridian $M^-(p)$ contains $o$ and has the same rank as $M$.

We assume from now on that $M$ is simply connected and irreducible. Then the isotropy group $K$ is connected. Of particular interest to us will be the so-called bottom space or adjoint space $\bar{M}$ of $M$. The bottom space $\bar{M}$ of $M$ is characterized by the property that every Riemannian symmetric space $M'$ that is locally isometric to $M$ is a Riemannian universal covering space of $\bar{M}$. It is constructed from $M$ by identifying all points with the same isotropy groups. The bottom space $\bar{M}$ has no poles and its geodesic symmetric are pairwise distinct isometries (see \cite [Lemma 2.1] {BOR19}). 

Let $\pi : M \to \bar{M}$ be the canonical projection. We put $\bar{p} = \pi(p) \in \bar{M}$ for $p \in M$ and $\bar{v} = d_p\pi(v) \in T_{\bar{p}}\bar{M}$ for $v \in T_pM$. The geodesic symmetry of $\bar{M}$ at $\bar{p}$ is denoted by $\bar{\sigma}_{\bar{p}}$. Let $\gamma_{\bar{v}} : [0,1] \to \bar{M}$ be a closed geodesic in $\bar{M}$ with period $1$ and $\gamma_{\bar{v}}(0) = \bar{o}$ $ (= \gamma_{\bar{v}}(1))$. Then $\bar{p} = \gamma_{\bar{v}}(\frac{1}{2})$ is an antipodal point of $\bar{o}$ and $\bar{M}^+(\bar{p})$ is a polar of $\bar{M}$ with $\dim(\bar{M}^+(\bar{p})) \geq 1$. The isometry $g^{\bar{v}} = \bar{\sigma}_{\bar{p}} \circ \bar{\sigma}_{\bar{o}} = \bar{\sigma}_{\bar{o}} \circ \bar{\sigma}_{\bar{p}} \in I(\bar{M})$ is involutive, nontrivial and fixes every point on $\gamma_{\bar{v}}([0,1])$. The linear isometry $\ell^{\bar{v}} = d_{\bar{o}}g^{\bar{v}}$ of $T_{\bar{o}}\bar{M}$ is involutive, nontrivial and coincides with parallel transport along $\gamma_{\bar{v}}$ from $\bar{o} = \gamma_{\bar{v}}(0)$ to $\bar{o} = \gamma_{\bar{v}}(1)$. In particular, $\ell^{\bar{v}}(\bar{v}) = \bar{v}$. Note that $g^{\bar{v}} = g^{\bar{w}}$ for any $\bar{w} \in T_{\bar{o}}\bar{M}$ with $\gamma_{\bar{w}}(\frac{1}{2}) = \bar{p}$.

The bottom space $\bar{M}$ can be written as $\bar{M} = G/\bar{K}$ with $\bar{K} = G_{\bar{o}}$. The identity component $\bar{K}^o$ of $\bar{K}$ is (isomorphic to) $K$. We now consider the isotropy representation of $\bar{K}^o \cong K$ on $T_{\bar{o}}\bar{M}$. For each $\bar{v} \in T_{\bar{o}}\bar{M}$, the linear isometry $\ell^{\bar{v}}$ leaves the isotropy orbit $K \cdot \bar{v} = K/K_{\bar{v}}$ invariant (see  \cite[Proposition 2.4]{BOR19} and its proof). For all $k \in K$ we have $g^{\overline{d_ok(v)}} = kg^{\bar{v}}k^{-1}$ and therefore $g^{\bar{v}} = kg^{\bar{v}}k^{-1}$ for all $k \in K_{\bar{v}}$. This implies $\ell^{\bar{v}} \circ d_{\bar{o}}k = d_{\bar{o}}k \circ \ell^{\bar{v}}$ for all $k \in K_{\bar{v}}$. The subspace $\Fix_{T_{\bar{o}}\bar{M}}(\ell^{\bar{v}})$ of fixed vectors of $\ell^{\bar v}$ in $T_{\bar{o}}\bar{M}$ always contains the normal space $\nu_{\bar v}(K\cdot {\bar v})$. Moreover, if $\rk(M) \geq 2$, then $\Fix_{T_{\bar{o}}\bar{M}}(\ell^{\bar{v}}) = \nu_{\bar v}(K\cdot {\bar v})$ if and only if $K \cdot {\bar v}$ is an extrinsically symmetric orbit.

For any Riemannian symmetric space $M' = G/K'$ that is locally isometric to $M$ we can make similar constructions, using the canonical projection $\pi' : M \to M'$ and defining $o', v'$ and other entities analogously. The subspace $\Fix_{T_{o'}M'}(\ell^{v'})$ of fixed vectors of $\ell^{v'}$ in $T_{o'}M'$ always contains the normal space $\nu_{v'}(K\cdot {v'})$. However, it may happen that $\Fix_{T_{o'}M'}(\ell^{v'}) =  T_{o'}M'$.

Let $\ca$ be a maximal abelian subspace of $\cp \cong T_oM$ and let $\alpha_1, \ldots,\alpha _r \in \ca^*$ be simple roots for the corresponding root system. Let $H^1,\ldots, H^r \in \ca$ be the dual basis of $\alpha_1,\ldots,\alpha_r$. We fix an index $i \in \{1,\ldots,r\}$. Note that $H^i$ belongs to a simplex of dimension $1$ in the closure $\bar C$ of the Weyl chamber $C$ that corresponds to the choice of the simple roots $\alpha_1,\ldots,\alpha_r$. In fact, $H^i$ belongs to any of the reflection hyperplanes associated with $\alpha _j$, $j\neq i$. Then the abelian part of the normal space $\nu_{H^i} (K \cdot H^i) = \cz_\cp(H^i)$ is the real span of $H^i$ and hence of dimension $1$. This means that $K\cdot H^i$ is a most singular (or focal) orbit of $K$. Conversely, if $K \cdot v$ is a most singular orbit, then $K \cdot v = K\cdot H^i$ for some $i\in \{1, \ldots , r\}$ with a suitable rescaling of $v$. In fact, by applying a suitable transformation in the Weyl group to $v$, we can assume that $v$ belongs to $\bar C$ and so $\alpha _j(v) \geq 0$ for all $j \in \{1, \ldots , r\}$. Then, since $K \cdot v$ is most singular, $ \alpha _j (v) = 0$ except for one index $i \in \{1, \ldots , r\}$. 
  
Since the real span of $H^i$ is the abelian part of $\cz_\cp(H^i)$, it is the tangent space at $o$ of a $1$-dimensional flat $S^1$ of $M$. After a suitable rescaling of $H^i$ to some vector $v$ we can assume that $\gamma _{v}: [0,1] \to M$ is a closed geodesic with period $1$. The same can be done if we replace $M$ by an arbitrary globally symmetric quotient $M'$ of $M$ and, in particular, by the bottom space $\bar{M}$. After replacing $\overline{H^i} = d_o\pi(H^i)$ by a scalar multiple $\bar{v}$, we obtain a closed geodesic $\gamma_{\bar{v}}: [0,1] \to \bar{M}$ of period $1$. 
  
Let $\alpha = \delta _1 H^1 + \ldots + \delta_r H^r$ be the highest root and assume that $\delta_i > 1$. Then the orbit $K \cdot H^i \cong  K \cdot \overline{H^i}$ is not an extrinsically symmetric orbit (see \cite[page 199]{BO16} or \cite{KN64}). Then the tangent space $T_{\bar o}\bar{M}^-(\bar{p})$ contains properly $\nu _{\bar v}(K \cdot \bar{v}) = \cz_\cp(\bar{v})$, where $\bar{p} = \gamma_{\bar{v}}(1/2)$.  Moreover, the meridian $\bar{M}^-(\bar{p})$ is semisimple. In fact, if it were not semisimple, choose $0 \neq \bar{w} \in T_{\bar{o}}\bar{M}^-(\bar{p})$ with $[\bar w, T_{\bar o}\bar{M}^-(\bar p)]=\{0\}$. Then $[\bar w, Z_\cp(\bar v)] = \{0\}$ and so $\bar w=\bar v$  up to rescaling, since the abelian part of $\cz_\cp(\bar v)$ is $1$-dimensional. Then 
\[
T_{\bar o}\bar{M}^-(\bar p) \subset \cz_\cp(\bar w) = \cz_\cp(\bar v)
\]
and therefore $T_{\bar o}\bar{M}^-(\bar p) = \cz_\cp(\bar v)$, which is a contradiction. 
    
We can replace $\bar M$ by any  globally symmetric quotient $M'$ of $M$ and obtain a similar result, but in this situation it may happen that $M'^-(p') = M'$ (namely, if $p'$ is a pole of $o'$).  So we have the following result that will be useful for our purposes, since totally geodesic submanifolds of simply connected symmetric spaces are not in general simply connected, but globally symmetric. 
    
\begin{prop} \label{ssm} 
Let $M=G/K$ an irreducible Riemannian symmetric space of compact type and $M' = G/K'$ be a symmetric quotient of $M$. Let $0 \neq v' \in T_{o'}M'$ be such that the isotropy orbit $(K')^o \cdot v'$ is a most singular and not extrinsically symmetric orbit. Then $\gamma_{v'}$ is a closed geodesic, which we may assume to be of period $1$.  Moreover, if $p'= \gamma_{v'}(1/2)$, then the tangent space $T_{o'}M'^-(p')$ of the meridian $M'^-(p')$ is a semisimple Lie triple system that contains properly the normal space $\nu_{v'}((K')^o\cdot v') = \cz_\cp(v)$.
\end{prop}
    
We now turn our attention to symmetric spaces whose root system is of type $A$. For $r \geq 2$, the irreducible simply connected Riemannian symmetric spaces of compact type whose root system is of type $A_r$ are $SU_{r+1}/SO_{r+1}$, $SU_{r+1}$, $SU_{2r+2}/Sp_{r+1}$ and $E_6/F_4$ (for which $r = 2$).

\begin{prop} \label{irr} 
Let $M=G/K$ be an irreducible Riemannian symmetric space with root system of type $A_r$, $r \geq 2$. Let $v \in T_oM \cong \cp$ so that $K \cdot v$ is a principal orbit in $\cp$, and hence isoparametric. Let $E(v)$ be the tangent space at $v$ of a curvature sphere $S(v)$ of $K \cdot v$. Then the connected isotropy group $(K_v)^o$ restricted to $E(v)$ acts irreducibly.
\end{prop}

\begin {proof} 
The key fact for the proof is that the subgroup of $K$ that acts on the curvature sphere $S(v)$, when restricted to $S(v)$, is the full group of isometries of $S(v)$. Since the isotropy representations for dual symmetric spaces are equivalent, we may assume that that $M$ is of  noncompact type. 
 
The normal space $\nu_v (K\cdot v)$ of $K \cdot v$ at $v$ is the maximal abelian subspace of $\cp$ that contains the regular tangent vector $v$. Moreover, the Weyl group $W$ associated with the isoparametric submanifold $K\cdot v$ coincides with the Weyl group of $M$ (corresponding to the root system determined by the maximal abelian subspace $\nu _v (K\cdot v)$). 

Let $\xi \in \nu _v(K\cdot v)$ so that $E(v)$ coincides with the $+1$-eigenspace of the shape operator $A_\xi$ of $K\cdot v$.  Equivalently, $u = v + \xi$ belongs to the reflection hyperplane $H$ of $\nu _v(K\cdot v)$ associated with $E(v)$, and does not belong to any other reflection hyperplane of the Weyl group. The tangent vector $u$ belongs to a simplex of dimension $r-1$ of the closure of a Weyl chamber (such a simplex is an open subset of $H$). Note that the curvature normal $\eta(v)$ associated with $E(v)$ is perpendicular to the hyperplane $H$ of $\nu _v (K\cdot v)$.
 
The focal parallel orbit $K\cdot u$ of $K \cdot v$ is a subprincipal orbit and
\[
\nu _u (K\cdot u) =\cz_\cp(u),
\]
where $\cz_\cp(u)= \{z\in \cp: [u, z]=0\}$ is the centralizer of $u$ in $\cp \cong T_oM$. The Lie triple system $\nu _u (K\cdot u)$ splits as
\[
\nu _u (K\cdot u) = H \oplus (E(v) \oplus \RR \eta (v)), 
\]
where $H$ is the abelian part of the Lie triple system and $E(v) \oplus \RR \eta (v)$ is a nonabelian Lie triple system of rank $1$. Note that $v \in \nu_u (K\cdot u)$, since $H \oplus \RR\eta (v) = \nu _v(K\cdot v)$, and so the rank of $\nu _u(K\cdot u)$ is $r$.  
 
Let $N \subset M$ be the symmetric space of rank $1$ associated with $E(v)\oplus \mathbb R \, \eta (v)$. The marked Dynkin diagram of $N$ consists of one of the nodes of the Dynkin diagram of $M$ with corresponding multiplicity $m$ (see e.g.\ \cite {HO92}). In our particular situation there are no double roots and  $m \in \{1,2,4,8\}$. Therefore, $N = G'/K'$ is an $(m+1)$-dimensional real hyperbolic space, where $G' \subset G$ are the glide transformations of $N$ (and so $K' \subset K$). If $X \subset M$ is the totally geodesic submanifold with $T_oX = \nu_u(K\cdot u)$, then 
\[
X = \RR^{r-1} \times N = (\RR^{r-1}\times G')/K'.
\]
 
The image of the representation of $(K_u)^o$  on $T_oX = \nu _u(K \cdot u)$ coincides with the image of the isotropy representation of $K'$ on $T_oX$ (see \cite [Theorem 2]{HO92}). Note that $K_v \subset K_u$, since $K \cdot v$ is a principal orbit and $u= v + \xi \in \nu _v(K\cdot v)$.
Hence $(K_u)_v = K_v$ and then 
\[
\left.(K_v)^o\right|_{\nu _u(K\cdot u)} 
= \left.((K_u)_v)^o\right|_{\nu _u(K\cdot u)}
= \left.(K'_v)^o\right|_{\nu _u(K\cdot u)},
\]
where we regard $K'$ as the isotropy group at $o$ of $X$, acting on $T_oX = H \oplus (E(v) \oplus \RR \eta (v)$. In particular, 
\begin {equation}\label {1}
\left.(K_v)^o\right|_{E(v)} 
= \left.(K'_v)^o\right|_{E(v)}.
 \end {equation}
Note that $\langle \eta (v), -v \rangle =1$ since $A_{-v}= \id$. We write $v = c\eta (v) + w$ with $c \neq 0$ and $w\in H$. Since $K'$ acts trivially on $H$, we get $K'_v = K' _{\eta (v)}$ and so $E(v)$ may be regarded as the tangent space of  the unit sphere of $T_oN$ at $\frac {1}{\Vert \eta (v)\Vert}\eta (v)$ . Since $N$ is an $(m+1)$-dimensional real hyperbolic space, we get $K' \cong SO_{m+1}$ and so $\left.K'_v\right|_{E(v)} \cong SO_m$, which acts irreducibly on $E(v)$. Then, from (\ref {1}), $(K_v)^o$ acts irreducibly on $E(v)$.
\end {proof}

The following proposition follows from the classification of polars (see \cite {NT00}, though the meridians are not explicitely listed there). Nevertheless, we include a direct proof here based on general arguments. Let $\Phi$ be a root system in $\RR^{r+1}$ of type $A_r$. We can assume that $\Phi = \{e_i - e_j :   i,j \in \{1, \ldots , r+1\}, \ i\neq j\}$, where $e_1, \ldots , e_{r+1}$ is the canonical basis of $\RR^{r+1}$.  Then $\Delta = \{ e_i - e_{i+1}: 1 \leq i \leq r\}$ is a set of simple roots of $\Phi$ and $\Phi ^+ = \{e_i - e_j :   1 \leq i < j\leq r +1\}$ is the resulting set of positive roots. Note that $e_1-e_2, \ldots , e_r-e_{r+1}$ is a basis of the  hyperplane $\RR_0^{r+1} = \{x \in \RR^{r+1} : x_1 + \ldots + x_{r+1} = 0\}$ of $\RR^{r+1}$. The Weyl group $W$ generated by the hyperplane reflections associated with $\Phi$ acts irreducibly on $\RR_0^{n+1}$. 
 
Let $\tilde{\Phi}$ be a nonempty proper root subsystem of $\Phi$. Then, up to a suitable relabelling, 
$\tilde{\Phi}$ is characterized in the following way: there exist integers $0 = d_0 < d_1< \ldots < d_k \leq r+1$ with $d_{\nu+1} - d_\nu \geq 2$ such that 
\[
 \tilde{\Phi} = \{e_i-e_j : i,j \in \{d_\nu +1, \ldots , d_{\nu+1}\}, \ i\neq j,\ 0 \leq \nu \leq k-1\}. 
\]
 The key fact for proving the above well-known equality is that if $e_i-e_j$ and $e_j -e_l$ belong to $\tilde{\Phi}$, then $e_i-e_l$ belongs to 
 $\tilde{\Phi}$ as well, which follows by applying to $e_j-e_l$ the hyperplane reflection of $\RR_0^{r+1}$ determined by $e_i-e_j$. Note that $\tilde{\Phi}$ is the direct sum of root systems of type $A_{d_{\nu+1}-(d_\nu +1)}$, $\nu =0, \ldots , k-1$.

Consider the following set ($I$) of $k + r -d_k + 1 \geq 2$ linearly independent linear equations of 
$\RR^{r+1}$:  
\[
(I) = \begin{cases}
x_{d_\nu +1} + x_{d_\nu + 2} + \ldots + x_{d_{\nu+1}} = 0 & ,\ \nu = 0, \ldots  , k-1, \\
x_{\mu+1} = 0 & ,\ \mu = d_k,\ldots,r.
\end {cases}
\]
If $V \subset \RR^{r+1}_0$ is the subspace determined by (I), then 
\[
\tilde{\Phi} = \Phi \cap V.
\]
Any $\alpha \in \Phi \setminus \tilde{\Phi}$ projects nontrivially onto $V^\perp \cap \RR^{r+1}_0$. Then, since $\Phi \setminus \tilde{\Phi}$ is finite, there exists $z\in V^\perp \cap \RR^{r+1}_0$ such that 
\begin {equation} \label{3} 
\tilde{\Phi} = \{\alpha \in \Phi : \langle \alpha , z \rangle = 0\}.
\end {equation}
We will use this equality in the proof of the following proposition.

\begin{prop} \label{A_r} 
Let $M = G/K$ be a simply connected Riemannian symmetric space of compact type with $r = \rk(M) \geq 2$ and root system of type $A_r$. Let $\bar M = G/\bar K$ be the bottom space of $M$. Then $V \subset T_{\bar o}\bar M$ is the tangent space to a meridian of $\bar M$ if and only if $V$ is the normal space to an extrinsically symmetric isotropy orbit of $\bar K^o \cong K$ (or equivalently, for our spaces, to a most singular orbit). In particular, $V$ is a nonsemisimple Lie triple system of  $T_{\bar o}\bar M$.
\end{prop}

\begin{proof}
Let $\bar{M}$ be the bottom space of $M$ and $\gamma_{\bar{v}}: [0,1] \to \bar{M}$ be a closed geodesic with period $1$ and $\gamma_{\bar{v}}(0)= \bar{o} = \gamma_{\bar{v}}(1)$. Consider the polar $\bar{M}^+(\bar{p})$ and the meridian $\bar{M}^-(\bar{p})$ through the antipodal point $\bar{p} = \gamma_{\bar{v}}(\frac{1}{2})$ of $\bar{o}$. The meridian $\bar{M}^-(\bar{p})$ is the connected component containing $\bar{o}$ of $\Fix_{\bar{M}}(g^{\bar{v}})$ and $T_{\bar{o}}\bar{M}^-(\bar{p}) = \Fix_{T_{\bar{o}}\bar{M}}(\ell^{\bar{v}})$. Recall that $\bar{M} = G/\bar{K}$ with $\bar{K} = G_{\bar{o}}$, $\bar{K}^o \cong K$, $g^{\bar{v}} = kg^{\bar{v}}k^{-1}$ and $\ell^{\bar{v}} \circ d_{\bar{o}}k = d_{\bar{o}}k \circ \ell^{\bar{v}}$ for all $k \in (\bar{K}^o)_{\bar{v}}$, and $\ell^{\bar{v}}(\bar{K} \cdot \bar{v}) = \bar{K} 
\cdot \bar{v}$.

We now choose $\bar{w} \in \bar{v} + \nu_{\bar{v}}(\bar{K}^o \cdot \bar{v})$ such that the orbit $\bar{K}^o \cdot \bar{w}$ is principal, and thus isoparametric. According to \cite[Proposition 2.4]{BOR19} and its proof, the linear isometry $\ell^{\bar v}$ fixes pointwise the normal space $ \nu_{\bar{v}}(\bar{K}^o \cdot \bar{v})$. In fact, the parallel transport $\ell^{\bar{v}}$ along $\gamma_{\bar{v}}$ must be trivial when restricted to any abelian subspace of $T_{\bar{o}}\bar{M}$ containing $\bar{v}$. We have $(\bar{K}^o)_{\bar{w}} \subseteq (\bar{K}^o)_{\bar{v}}$ since $\bar{K}^o \cdot \bar{w}$ is a principal orbit. Moreover, we have $\nu_{\bar{w}}(\bar{K}^o \cdot \bar{w}) \subseteq \nu_{\bar{v}}(\bar{K}^o \cdot \bar{v})$ by construction. Then $\ell^{\bar{v}}(\bar{K}^o \cdot \bar{w}) = \bar{K}^o \cdot \bar{w}$, $\ell^{\bar{v}} \circ d_{\bar{o}}k = d_{\bar{o}}k \circ \ell^{\bar{v}}$ for all $k \in (\bar{K}^o)_{\bar{w}}$, and $\ell^{\bar{v}}$ is the identity on $\nu_{\bar{w}}(\bar{K}^o \cdot \bar{w})$.

Let $S$ be the connected component of $\Fix_{T_{\bar{o}}\bar{M}}(\ell^{\bar{v}}) \cap (\bar{K}^o \cdot \bar{w})$ containing $\bar{w}$. Since $\ell^{\bar{v}}$ is an isometry of the orbit $\bar{K}^o\cdot \bar{w}$, we see that $S$ is a totally geodesic submanifold of $\bar{K} ^o \cdot \bar{w}$. Let $W$ be the orthogonal complement of $T_{\bar{w}}S$ in $T_{\bar{w}}(\bar{K}^o \cdot \bar{w})$ and $\alpha$ be the second fundamental form of $\bar{K}^o \cdot \bar{w}$. For $X \in W$ and $Y \in T_{\bar{w}}S$ we have
\[
\alpha (X,Y) =  \ell^{\bar{v}} (\alpha(X,Y)) = \alpha (\ell^{\bar{v}}(X),\ell^{\bar{v}}(Y)) = \alpha(-X,Y) = -\alpha (X,Y),
\]
and hence $\alpha (X,Y)=0$. It follows that $T_{\bar{w}}S$ is invariant under the shape operator $A_{\bar{\xi}}$ of $\bar{K}^o \cdot \bar{w}$ for any normal vector $\bar{\xi}$ of $S$ at $\bar{w}$. This implies that $S$ is an isoparametric submanifold of $T_{\bar{o}}\bar{M}$ (which is contained in $\bar{w} + W^\perp$). 

If $E'(\bar{w})$ is the tangent space at $\bar{w}$ of a curvature sphere of $S$, then $E'(\bar{w})$ is contained in the tangent space $E(\bar{w})$ of some curvature sphere of $\bar{K}^o \cdot \bar{w}$. Since $S$ is invariant under the action of $(\bar{K}^o)_{\bar{w}}$, also $E'(\bar{w})$ is invariant by $(\bar{K}^o)_{\bar{w}}$. From Proposition \ref{irr} we know that $(\bar{K}^o)_{\bar{w}}$ acts irreducibly on  $E(\bar{w})$, and consequently
\[
E'(\bar{w}) = E(\bar{w}). 
\]
Therefore, any curvature sphere of $S$ is a curvature sphere of $\bar{K}^o \cdot \bar{w}$. This implies that the Weyl group of $S$ is generated by some of the reflection hyperplanes of the Weyl group associated with $\bar{K}^o \cdot \bar{w}$. Note that the set  $J' = \{\eta'_1(\bar{w}),\ldots,\eta'_d(\bar{w})\}$ of the curvature normals $S$ at $\bar{w}$ is a subset of the set $J$ of the curvature normals of $\bar{K}^o \cdot \bar{w}$ at $\bar{w}$. 

From (\ref {3}) and its preceding paragraph, there exists $\bar{z} \in \nu_{\bar{w}}(\bar{K}^o \cdot \bar{w})$ such that $\eta \in J$ satisfies $\langle \eta , \bar{z} \rangle = 0$ if and only if $\eta \in J'$. This implies that $T_{\bar{w}}S$ coincides with the $(+1)$-eigenspace of the shape operator of $\bar{K}^o \cdot \bar{w}$ with respect to $-\bar{w}+\bar{z}$. Note that $\bar{z} = \bar{w} + (-\bar{w} +\bar{z})$ and so $\bar{K}^o \cdot \bar{z}$ is a parallel focal orbit to $\bar{K}^o \cdot \bar{w}$ and $S$ is the connected component of the fibers of the parallel map from $\bar{K}^o \cdot \bar{w}$ into $\bar{K}^o \cdot \bar{z}$. Then 
\[
\nu_{\bar{z}}(\bar{K}^o \cdot \bar{z}) = T_{\bar{w}}S \oplus \nu_{\bar{w}}(\bar{K}^o \cdot \bar{w}). 
\]

Note that $\ell^{\bar{v}}(\bar{z}) = \bar{z}$ and $\ell^{\bar{v}}(\bar{K}^o \cdot \bar{z}) = \bar{K}^o \cdot \bar{z}$. Moreover, $\ell^{\bar{v}}$ is the identity on $\nu_{\bar{z}}(\bar{K}^o \cdot \bar{z})$ and minus the identity on $T_{\bar{z}}(\bar{K}^o \cdot \bar{z})$. Thus we see that $\bar{K}^o \cdot \bar{z}$ is extrinsically symmetric. Since the set of fixed vectors of $\ell^{\bar{v}}$ coincides with $\nu_{\bar{z}} (\bar{K}^o \cdot \bar{z})$, the assertion follows.
\end{proof}

The next two results are not related to the above, but will be useful later.

\begin{prop} \label{indexSigmaM}
Let $M$ be an irreducible Riemannian symmetric space and $\Sigma$ be a connected totally geodesic submanifold of $M$ with $\dim(\Sigma) < \dim(M)$. Then $i(\Sigma) \leq i(M)$.
\end{prop}

\begin{proof}
Let $\Sigma'$ be a maximal totally geodesic submanifold of $M$ with $d = \codim(\Sigma') = i(M)$. We can assume that $o \in \Sigma \cap \Sigma'$. By assumption, the isotropy group $K$ of $M = G/K$ acts irreducibly on $T_oM$. Therefore, for all $0 \neq X \in T_oM$ the span of $\{d_ok(X) : k \in K\}$ is equal to $T_oM$. It follows that there exists $k \in K$ such that $T_o\Sigma$ is not contained in $d_ok(T_o\Sigma')$, or equivalently, $\Sigma$ is not contained in $k(\Sigma')$. Then $\Sigma \cap k(\Sigma')$ is a proper totally geodesic submanifold of $\Sigma$ and the codimension of $\Sigma \cap k(\Sigma')$ in $\Sigma$ is less than or equal to $d$. Consequently, $i(\Sigma) \leq d = i(M)$.
\end{proof}

\begin{re} \label{obstruction}
{\rm Let $M = G/K$ and $\Sigma = G'/K'$ be irreducible Riemannian symmetric spaces of noncompact type. We have the following obvious necessary conditions for the existence of a totally geodesic embedding of $\Sigma$ into $M$:
\begin{align*}
& \dim(\Sigma) < \dim(M),\ \dim(G') < \dim(G),\ \dim(K') < \dim(K), \\
& \rk(\Sigma) \leq \rk(M),\ \rk(G') \leq \rk(G),\ \rk(K') \leq \rk(K).
\end{align*}
Proposition \ref{indexSigmaM} shows that the index imposes further necessary conditions:
\[
i(\Sigma) \leq i(M) \leq \dim(M) - \dim(\Sigma).
\]
Choose for example $M = SO^o_{d, d +k}/SO_d SO_{d+k}$. For given $\Sigma$, if we choose $k$ sufficiently large and $d \geq \rk(\Sigma)$, then the first set of necessary conditions is satisfied. We know from \cite{BO18} that $\rk(\Sigma) \leq i(\Sigma)$, and from \cite{BO16} that equality holds only for $SL_{r+1}(\RR)/SO_{r+1}$ and $SO^o_{r, r +l}/SO_r SO_{r+l}$. If $\Sigma$ is different from these symmetric spaces, then we can choose $d$ with $\rk(\Sigma) \leq d < i(\Sigma)$. Then $i(\Sigma) > d = i(M)$ and therefore the second set of necessary conditions tells us that there cannot be a totally geodesic embedding of $\Sigma$ into $SO^o_{d, d +k}/SO_d SO_{d+k}$ for any $k \geq 0$. Thus the index gives a useful additional obstruction for the existence of totally geodesic embeddings (or immersions in the compact case) in addition to the standard obstructions given by dimensions and ranks.} 
\end{re}

\begin{re}
{\rm It is known that every Riemannian symmetric space $M$ of noncompact type admits a totally geodesic embedding into $SL_n(\RR)/SO_n$ for some $n \in \NN$. This is a consequence of the well-known unitary trick (and can also be seen as a particular case of Karpelevich's Theorem by embedding $I(M)^o$ into the special linear group via the adjoint representation). By duality, any symmetric space $M$ of compact type admits a totally geodesic immersion into $SU_n/SO_n$ for some $n \in \NN$. The symmetric space $SU_n/SO_n$ is a totally geodesic hypersurface of the symmetric space $U_n/SO_n$. The symmetric space $U_n/SO_n$ is a symmetric $R$-space, arising as a symmetric orbit of the isotropy representation of the symmetric space $Sp_n/U_n$. Hence $U_n/SO_n$ admits an isometric immersion into $\RR^{n(n+1)}$ with parallel second fundamental form. Then, by a well-known result of Vilms \cite{V72}, the corresponding Gauss map is totally geodesic. Consequently, every Riemannian symmetric space of compact type admits a totally geodesic immersion into some real Grassmannian $SO_{2r+k}/SO_rSO_{r+k}$. Then $i(M) \leq i (SO_{2r+k}/SO_rSO_{r+k}) = r$. Thus the index of $M$ is a lower bound for the rank $r$ of a real Grassmannian into which $M$ can be totally geodesically immersed (where the index of $M$ is defined to be the sum of the indices of the 
locally irreducible components of $M$).}
\end{re} 

The next result is useful for the investigation of totally geodesic submanifolds in reducible Riemannian symmetric spaces.

\begin{prop} \label{productestimate}
Let $M = M_1 \times M_2$ be the Riemannian product of two irreducible Riemannian symmetric spaces $M_1,M_2$ 
and $\Sigma = G'/K'$ be a totally geodesic submanifold of $M$. Let $o = (o_1,o_2) \in \Sigma$ and assume that $T_{o_1}M_1 \times \{0\}$ and $\{0\} \times T_{o_2}M_2$ are not contained in $T_o\Sigma \subseteq T_oM = T_{o_1}M_1 \times T_{o_2}M_2$. Then $\codim(\Sigma) \geq i(M_1) + i(M_2)$.
\end{prop}

\begin{proof}
Let $\pi_j : M = M_1 \times M_2 \to M_j$ be the canonical projection. 

If $\pi _1(\Sigma) \subsetneq M_1$ and $\pi_2(\Sigma) \subsetneq M_2$, then $\Sigma$ is contained in the totally geodesic submanifold $\pi_1(\Sigma) \times \pi _2(\Sigma)$ of $M$ and 
\begin{align*}
\codim_M(\Sigma) & \geq  \codim_M(\pi_1(\Sigma) \times \pi _2(\Sigma))  = \codim_{M_1}(\pi_1(\Sigma)) + \codim_{M_2}(\pi_2(\Sigma)) \\
& \geq i (M_1) + i(M_2).
\end{align*} 
We can therefore assume, without loss of generality, that $\pi _1(\Sigma) = M_1$. We define two subspaces $V_j = \ker(d_o\pi_j|_{T_o\Sigma})$ ($j =1,2$) of $T_o\Sigma$ and note that both are $K'$-invariant.

{\sc Case 1: $V_1 = \{0\}$.} 
Then $\pi_1|_\Sigma : \Sigma \to M_1$ is a totally geodesic local diffeomorphism and hence an affine map. It follows that $\pi_1|_\Sigma$ is a homothety and therefore $\Sigma$ is an irreducible symmetric space. 

If  $\dim (M_1) \leq \dim (M_2)$, then $\dim(\Sigma) = \dim (M_1) \leq \frac{1}{2}\dim (M)$ and
\[
\codim(\Sigma) \geq \textstyle{\frac{1}{2}}\dim(M) = \textstyle{\frac{1}{2}}\dim(M_1) + \textstyle{\frac{1}{2}}\dim(M_2) \geq  i(M_1) + i(M_2), 
\]
where the last inequality follows from the fact that every irreducible symmetric space contains at least one pair of perpendicular reflective submanifolds (\cite{L75}, \cite{L79}).

If $\dim(M_1) > \dim(M_2)$, then $V_2 \neq \{0\}$ because of $\dim(\Sigma) = \dim(M_1)$. Since $V _2$ is $K'$-invariant and $\Sigma$ is irreducible, it follows that $V_2 = T_o\Sigma$. Then $\Sigma = M_1$ and hence $T_o\Sigma = T_{o_1}M_1 \times \{0\}$, which contradicts the assumption. 

{\sc Case 2: $V_1 \neq \{0\}$.}
Then $\Sigma$ is a Riemannian product $\Sigma = \Sigma _1 \times \Sigma _2$, where $\Sigma_1$ is homothetic to $M_1$ and $\pi_2|_{\Sigma_2} : \Sigma _2 \to M_2$ is a totally geodesic immersion with $T_o\Sigma _2 = V_1$. If $\pi_2|_{\Sigma_1} : \Sigma_1 \to M_2$ is constant, then $\Sigma_1 = M_1$, which contradicts the assumption. Let us consider the totally geodesic map $\pi_2|_{\Sigma} : \Sigma \to M_2$. Assume that  
$V _2 \neq \{0\}$. Since $(\pi_1,\pi _2)|_\Sigma$ is the inclusion map of $\Sigma$ into $M$, we have $V_1 \cap V_2=\{0\}$. Then, as $V _2$ is $K'$-invariant, we must have $V_2= T_o\Sigma _1$, which implies $T_o\Sigma _1 = T_{o_1}M_1$ and contradicts the assumption. Therefore $\pi_2|_\Sigma : \Sigma \to M_2$ is a totally geodesic immersion. If $\pi_2(\Sigma) = M_2$, then $\pi_2|_\Sigma : \Sigma \to M_2$ is an affine local diffeomorphism and so $M_2 =  \Sigma = \Sigma_1\times \Sigma_2$ up to rescaling of the metric in each irreducible factor. This is a contradiction since $M_2$ is irreducible. Thus $\pi_2(\Sigma)$ is strictly contained in $M_2$ and so \[
\dim(\Sigma) \leq \dim (M_2) - i(M_2) \leq \dim (M_1) - i(M_1) + \dim(M_2) - i(M_2),
\]
which implies $\codim(\Sigma) \geq i(M_1) + i (M_2)$.
\end{proof}

\section{Lagrangian Grassmannians} \label{LG}

The complex $2$-plane Grassmannian $G_2(\CC^{2r+2}) = SU_{2r+2}/S(U_2U_{2r})$ is the complexification of the quaternionic projective space $\HH P^r = Sp_{r+1}/Sp_1Sp_r$, or equivalently, $\HH P^r$ is a real form of $G_2(\CC^{2r+2})$. The symmetric space $M = SU_{2r+2}/Sp_{r+1}$ is the Lagrangian Grassmannian of all real forms of $G_2(\CC^{2r+2})$ that are congruent to $\HH P^r$. We have $\rk(M) = r$ and $\dim(M) = r(2r+3)$. The associated root system is of type $A_r$ and all roots have multiplicity $4$. The symmetric space $M = SU_4/Sp_2 \cong Spin_6/Spin_5$ is isometric to $S^5$ and hence $i(M) = 1$. Oni\v{s}\v{c}ik (\cite{O80}) proved that $i(M) = 6$ for $M = SU_6/Sp_3$. In this  section we will prove that $i(SU_8/Sp_4) = 11$ and $i(SU_{2r+2}/Sp_{r+1}) = 4r$ for $r \geq 4$. Throughout this section we assume that $r \geq 3$. We know from \cite{L79} that $i_r(SU_8/Sp_4) = 11$ and $i_r(SU_{2r+2}/Sp_{r+1}) = 4r$ for $r \geq 4$.

\begin{lm} \label {SigmaA_r}
Let $\Sigma$ be a maximal, locally reducible, totally geodesic submanifold of $M = SU_{2r+2}/Sp_{r+1}$. Assume that $\Sigma$ has a local de Rham factor $\Sigma^1 = G^1/K^1 \subset M$ whose root system is not of type $A_s$, where $s = \rk(\Sigma^1)$. Then $\Sigma$ is nonsemisimple. 
\end{lm}

\begin {proof}
Let $\bar{M}$ be the bottom space of $M$ and $\pi : M \to \bar{M}$ be the canonical projection. We use the notations introduced in Section \ref{gsr}. Let $\bar{\Sigma}$ be the maximal totally geodesic submanifold of $\bar{M}$ with $T_{\bar{o}}\bar{\Sigma} = d_o\pi(T_o\Sigma)$. We denote by $\bar{\Sigma}^1 = G^1/\bar{K}^1 \subset \bar{M}$ the corresponding local de Rham factor of $\bar{\Sigma}$ whose root system is not of type $A_s$. Locally, around $\bar{o}$, we can write $\bar{\Sigma} = \bar{\Sigma}^1 \times \bar{\Sigma}'$.

Assume that $\Sigma$ is semisimple. Then $\bar{\Sigma}$ is semisimple. Let  $\alpha_1, \ldots, \alpha_s$ be a set of simple roots of the root system of $\bar{\Sigma}^1$ and $\bar{H}^1,\ldots,\bar{H}^s \in T_{\bar{o}}\bar{\Sigma}^1$ be the dual basis of $\alpha_1, \ldots, \alpha_s$.  Let $\alpha = \delta _1\alpha _1 + \ldots + \delta_s\alpha _s$ be the highest root. Since the root system of $\bar{\Sigma}^1$ is not of type $A_s$, we have $\delta_i > 1$ for some $i \in \{1, \ldots , s \}$. Then $(\bar{K}^1)^o \cdot \bar{H}^i$ is a most singular and not extrinsically symmetric orbit in $T_{\bar{o}}\bar{\Sigma}^1$. 

Since $\bar{H}^i$ is a most singular vector, we can rescale $\bar{H}^i$ to a vector $\bar{v} \in T_{\bar{o}}\bar{\Sigma}^1$ so that the closed geodesic $\gamma_{\bar{v}} : [0,1] \to \bar{\Sigma}^1$ has period $1$ (see e.g.\ proof of \cite[Proposition 2.4]{BOR19}). Let $\bar{p} = \gamma_{\bar{v}}(\frac{1}{2})$ be the antipodal point of $\bar{o}$ on $\gamma_{\bar{v}}$. It follows from Proposition \ref{ssm} that the tangent space $T_{\bar{o}}(\bar{\Sigma}^1)^-(\bar{p})$ of the meridian $(\bar{\Sigma}^1)^-(\bar{p})$ is a semisimple Lie triple system. It may happen that $T_{\bar{o}}(\bar{\Sigma}^1)^-(\bar{p}) = T_{\bar{o}}\bar{\Sigma}^1$ if $\bar{\Sigma}^1$ has poles.

The meridian $(\bar{\Sigma}^1)^-(\bar{p})$ of $\bar\Sigma^1$ is contained in the meridian $\bar{M}^-(\bar{p})$ of $\bar{M}$. By Proposition \ref{A_r}, $V = T_{\bar{o}}\bar{M}^-(\bar{p})$ is a nonsemisimple Lie triple system of $T_{\bar{o}}\bar{M}$. Moreover, as explained in Section \ref{gsr}, $V$ contains the centralizer $\cz_{T_{\bar{o}}\bar{M}}(\bar{v})$ of $\bar{v}$ in $T_{\bar{o}}\bar{M}$. Note that $T_{\bar{o}}\bar{\Sigma}' \subset \cz_{T_{\bar{o}}\bar{M}}(\bar{v})$.

Note that 
\[
W  = T_{\bar{o}}(\bar{\Sigma}^1)^-(\bar{p}) \oplus T_{\bar{o}}\bar{\Sigma}' 
\]
is a semisimple Lie triple system in $T_{\bar{o}}\bar{M}$ which is contained in the nonsemisimple Lie triple system $V$. This implies that there exists $\bar{z} \in V$ with $\bar{z} \notin W$ so that $[\bar{z}, W] = \{0\}$. In particular,  $[\bar{z},  T_{\bar{o}}\bar{\Sigma}'] = \{0\}$. Then 
\[
U  = \cz_{T_{\bar{o}}\bar{M}}(T_{\bar{o}}\bar{\Sigma}') + T_{\bar{o}}\bar{\Sigma}'
\]
is a proper Lie triple system in $T_{\bar{o}}\bar{M}$. Since $\bar{z} \in U$, we see that $T_{\bar{o}}\bar{\Sigma}$ is properly contained in $U$. This is a contradiction to the maximality of $\bar{\Sigma}$.  It follows that $\Sigma$ is nonsemisimple. 
\end{proof}

\begin {cor} \label {corSigmaA_r}
Let $M = SU_{2r +2}/Sp_{r+1}$ and $\Sigma$ be a maximal totally geodesic submanifold of $M$. Assume that $\Sigma$ has a local de Rham factor whose root system is not of type $A$. Then $\Sigma$ is locally irreducible. 
\end{cor}

\begin {proof}
Assume that $\Sigma$ is locally reducible. Then $\Sigma$ is nonsemisimple by Lemma \ref {SigmaA_r}. From \cite[Theorem 1.2]{BO16} we see that $T_o\Sigma$ coincides with the normal space of an extrinsically symmetric isotropy orbit. However, from \cite[Table 3]{BO16} we know that the root system of any irreducible factor of a maximal nonsemisimple totally geodesic submanifold of $M$ is of type $A$. It follows that $\Sigma$ must be locally irreducible. 
\end {proof}

From the classification of Riemannian symmetric spaces and their root systems with multiplicities (see e.g.\ \cite[Section 13.1]{BCO16}) we immediately get the following result. 

\begin {lm}\label{AAA} 
Let $\Sigma = G'/K'$ be an irreducible simply connected Riemannian symmetric space of compact type whose root system is not of type $A_s$, $s = \rk(\Sigma)$. Let $m_1$ and $m_2$ be the associated multiplicities of $\Sigma$. Then one of the following statements holds: 
\begin{itemize}
\item[({\rm i})] $\Sigma = G'/K'$ is inner, that is, $\rk(G') = \rk(K')$;
\item [({\rm ii})] $\max\{m_1, m_2\}\leq 2$;
\item [({\rm iii})] $\Sigma = SO_{2s+n}/SO_sSO_{s+n}$ with $s \geq 3$ odd and $n \geq 4$ even. {\rm [}In this case $\Sigma$ is outer, $\rk(\Sigma) = s$ and $(m_1,m_2) = (1,n)$.{\rm ]}
\end{itemize}
\end{lm}

We will now investigate these three possibilities in more detail. 

\begin{prop} \label{grassmannian}
Let $\Sigma$ be a totally geodesic submanifold of $M = SU_{2r +2}/Sp_{r+1}$ and assume that $\Sigma$ is locally isometric to $SO_{2s+n}/SO_sSO_{s+n}$ with $s \geq 3$ odd and $n \geq 4$ even. Then $\codim(\Sigma) > i_r(M)$.
\end{prop}
 
\begin {proof}
From \cite[Table 5]{BO16} we know that
\[
\Sigma' = SO_{2s+n-1}/SO_sSO_{s+n-1}
\]
is a maximal totally geodesic submanifold of $SO_{2s+n}/SO_sSO_{s+n}$ for which the codimension is equal to the index of $SO_{2s+n}/SO_sSO_{s+n}$. Furthermore, $\Sigma'$ is a reflective submanifold and its complementary reflective submanifold $(\Sigma')^\perp$ is locally isometric to an $s$-dimensional sphere $S^s$. Note that $\Sigma'$ is an inner symmetric space because $s$ is odd and $n$ is even, and hence its geodesic symmetry $\tau$ at $o$ is in $SO_sSO_{s+n-1}$. Moreover, we have $d_o\tau = I_s \otimes (-I_{s+n-1})$. Note that $\tau$ is in the center $\{(I_s,\pm I_{s+n-1})\}$ of $SO_sSO_{s+n-1}$, because $s \geq 3$ is odd and $n \geq 4$ is even. Therefore $\tau$ can be considered as an involutive isometry of $SO_{2s+n}/SO_sSO_{s+n}$ with $d_o\tau(X) = X$ for all $X \in T_o(\Sigma')^\perp$. Geometrically, $\tau$ is the isometric reflection of $SO_{2s+n}/SO_sSO_{s+n}$ in the reflective submanifold $(\Sigma')^\perp$. 

We now consider this setup in the bottom space $\bar{M}$ of $M$ via the canonical projection $\pi : M \to \bar{M}$ with the corresponding totally geodesic submanifolds $\bar{\Sigma}$, $\bar{\Sigma}'$, $(\bar{\Sigma}')^\perp$ and involution $\bar{\tau}$. If $\bar{\Sigma}$ is a reflective submanifold of $\bar{M}$, then $\Sigma$ is a reflective submanifold of $M$ and the assertion is obvious. So let as assume that $\bar{\Sigma}$ is not a reflective submanifold of $\bar{M}$. Since $\bar{\tau}$ is an involutive isometry in the center of the identity component $(\bar{K}')^o$ of the isotropy group $\bar{K}'$ of $\bar{\Sigma}$, we can consider $\bar{\tau}$ as an involutive isometry of $\bar{M}$.

Let $V^+$ and $V^- $ be the $(+1)$- and $(-1)$-eigenspaces of $d_{\bar{o}}\bar{\tau}$, respectively. Note that $V^+$ and $V^- $ are complementary reflective Lie triple systems. Moreover, we have $T_{\bar{o}}\bar{\Sigma}' \subset V^-$ and $T_{\bar{o}}(\bar{\Sigma}')^\perp \subset V^+$ by construction of $\bar{\tau}$. If $T_{\bar{o}}\bar{\Sigma}'= V^- $, then $\bar{\Sigma}'$ is reflective and hence $\bar{\Sigma}$ is reflective by \cite[Corollary 2.9]{BOR19} (since $\bar{\Sigma}$ contains $\bar{\Sigma}'$), which contradicts the assumption that $\bar{\Sigma}$ is not a reflective submanifold of $\bar{M}$. Thus we have a nontrivial orthogonal decomposition 
\[
V^- = T_{\bar{o}}\bar{\Sigma}' \oplus V^-_1.
\]
If $\dim(V^-_1) \geq s$, then $V^-$ is a reflective Lie triple system with $\dim(V^-) \geq \dim(\bar{\Sigma})$. The reflective submanifold $\tilde{\Sigma}$ of $M$ with $d_o\pi(T_o\tilde{\Sigma}) = V^-$ then satisfies $\dim(\tilde{\Sigma}) \geq \dim(\bar{\Sigma}) = \dim(\Sigma)$. If $\dim (V^-_1) < s$, then the isotropy group $SO_s SO_{s+n-1}$ acts trivially on 
 $V^-_1$. Since $\bar{\tau}$ belongs to this isotropy group by construction, $d_{\bar{o}}\bar{\tau}$ is the identity on $V^-_1$, which is a contradiction to $V^-_1 \subset V^-$.  
  
Since $n\geq 4$, $SO_{s+n-1}$ (which contains the symmetry $\tau$) acts trivially on $\RR^d$ for all $d \leq s+2$. This implies that $\dim(\tilde{\Sigma}) \geq \dim(\Sigma) + 2 > \dim (\Sigma)$, and so $\codim(\Sigma) > i_r(M)$.
\end{proof}

\begin{prop} \label{mleq2} 
Let $\Sigma = G'/K'$ be a locally irreducible totally geodesic submanifold of $M= SU_{2r +2}/Sp_{r+1}$ ($r \geq 3$) with associated multiplicities $m_1$ and $m_2$. If $\max \{m_1, m_2\} \leq 2$, then $\codim(\Sigma) \geq 4r \geq i_r(M)$.
\end{prop}

\begin{proof}
Recall that the root system associated with $M$ is of type $A_r$ and all roots have multiplicity $4$.
The number of reflection hyperplanes of a Weyl group of type $A_r$ is $\frac {r(r+1)}{2}$ and coincides with the number $\vert \Phi ^+\vert$ of positive roots in the corresponding root system. 
Since every root has multiplicity $4$, we get $\dim(M) = 4\vert \Phi ^+\vert +r = 2r^2 +3r$. 

Let $l$ be the number of reflections hyperplanes of the Weyl group of the universal covering space $\tilde{\Sigma}$ of $\Sigma$. Then $l = \bar{l} + \bar{s}$, where $\bar{l}$ is the number of long positive roots and $\bar{s}$ is the number of short positive roots in the root system associated with $\tilde{\Sigma}$, taking into account our conventions made near the beginning of Section \ref{gsr}.
Then, using the assumption that $\max \{m_1, m_2\} \leq 2$ and the fact that $\rk(\Sigma) \leq \rk(M)$, we get
\[
\dim (\Sigma) \leq l(\max \{m_1, m_2\}) + \rk(\Sigma) \leq 2l + r .
\]
From Proposition \ref{cardinal} we know that $l \leq  \frac {r(r+1)}{2}$ and hence
\[
\dim (\Sigma) \leq 2\frac {r(r+1)}{2} + r = r^2 + 2r.
\]
This implies 
\[
\codim(\Sigma) = \dim(M) - \dim(\Sigma) \geq 2r^2 + 3r - r^2 - 2r = r^2 + r = r(r+1) \geq 4r,
\]
since $r\geq 3$. 
\end{proof}

\begin{prop} \label{factorsAAA}
Let $\Sigma$ be a maximal semisimple totally geodesic submanifold of $M= SU_{2r +2}/Sp_{r+1}$. If the root system of every local de Rham factor of $\Sigma$ is of type $A$, then $\codim(\Sigma) > 4r$. 
\end{prop}

\begin{proof}
We denote by $\tilde{\Sigma}$ the Riemannian universal covering space of $\Sigma$ and consider its de Rham decomposition $\tilde{\Sigma} = \tilde{\Sigma}_1 \times \ldots \times \tilde{\Sigma_b}$, $b \geq 1$. 

{\sc Case 1.} Assume that $\tilde{\Sigma}_1 = S^{k_1}$ for some $k_1 \geq 2$ (after a suitable relabelling of the factors). 

If $\tilde{\Sigma}$ is reducible, it follows from (the dual version of) Proposition \ref{hyperbolicfactors} that $\tilde{\Sigma} = S^{k_1} \times S^{k_2}$ with $k_2 \geq 2$. The rank of the isotropy group of $\tilde{\Sigma}$ must satisfy
\[
\rk(SO_{k_1} \times SO_{k_2}) \leq  \rk(Sp_{r+1}) = r+1.
\]
Since 
\[
\frac{1}{2} \dim(\Sigma) -1 = \frac{k_1 -1}{2} + \frac{k_2 -1}{2} \leq \rk(SO_{k_1}\times SO_{k_2}), 
\] 
we obtain 
\[
\dim (\Sigma)\leq  2\rk(SO_{k_1}\times SO_{k_2}) + 2 \leq 2r +4.
\]
From this we get
\[
\codim(\Sigma) = \dim(M) - \dim (\Sigma) \geq (2r^2 + 3r) - (2r+4) = 2r^2 + r - 4 > 4r
\]
since $r \geq 3$, which is a contradiction. 

If $\tilde{\Sigma}$ is irreducible, we have $\tilde{\Sigma} = S^{k_1}$. We must have $\left[\frac{k_1}{2}\right] = \rk(SO_{k_1}) \leq  \rk(Sp_{r+1}) = r + 1$, which gives $k_1 \leq 2r+3$ and hence
\[
\codim(\Sigma) = r(2r+3) - k_1 \geq r(2r+3) - (2r+3) = 2r^2 + r - 3 > 4r. 
\]
It follows that none of the de Rham factors of $\tilde{\Sigma}$ is a sphere.

{\sc Case 2.} Assume that $\tilde{\Sigma}_1 = E_6/F_4$ (after a suitable relabelling of the factors). Recall that $\rk(E_6/F_4) = 2$, $\dim(E_6/F_4) = 26$, and the associated multiplicities of $E_6/F_4$ are $m_1 = 8 = m_2$. It follows from Proposition \ref{multcomp} that $\rk(\tilde{\Sigma}) < \rk(M) = r$. 

If $r = 3$, then $\tilde{\Sigma} = E_6/F_4$, $\dim(\Sigma) = 26$ and $\dim(SU_8/Sp_4) = 27$. Since $SU_8/Sp_4$ does not admit a totally geodesic hypersurface, this case cannot occur.

If $r = 4$, then $\tilde{\Sigma} = E_6/F_4$ since $\tilde{\Sigma}$ cannot have a rank $1$ factor by Case 1. Then $\dim(\Sigma) = 26$, $\dim(SU_{10}/Sp_5) = 44$, and hence $\codim(\Sigma) = 18 > 16 = 4r$.

If $r \geq 5$, then $\tilde{\Sigma} = E_6/F_4 \times \tilde{\Sigma}'$, where $\tilde{\Sigma}' = \{0\}$ or  $\tilde{\Sigma}' = \tilde{\Sigma}_2 \times \ldots \times \tilde{\Sigma_b}$ and each $\tilde{\Sigma}_i$ is an irreducible, simply connected, Riemannian symmetric space with $\rk(\tilde{\Sigma}_i) \geq 2$ and root system of type $A$. The isotropy group $\tilde{K}$ of $\tilde{\Sigma}$ must satisfy $\rk(\tilde{K}) \leq \rk(Sp_{r+1}) = r+1$, and thus the isotropy group $\tilde{K}'$ of $\tilde{\Sigma}'$ must satisfy $\rk(\tilde{K}') \leq r-3$. From the list of symmetric spaces with root system of type $A$ we can easily find the symmetric spaces $\tilde{\Sigma}'$ of maximal possible dimension with $\rk(\tilde{\Sigma}') \leq r-3$ and $\rk(\tilde{K}') \leq r-3$. They are:
\begin{itemize}
\item[(i)] $\tilde{\Sigma}' = SU_3$ if $r = 5$;
\item[(ii)] $\tilde{\Sigma}' = SU_4$ if $r = 6$;
\item[(iii)] $\tilde{\Sigma}' = SU_{2(r-3)}/Sp_{r-3}$ if $r \geq 7$.
\end{itemize}
In particular, such a $\tilde{\Sigma}'$ is always irreducible.

If $r = 5$, then $\dim(\Sigma) \leq \dim(E_6/F_4) + \dim(SU_3) = 34$. Since $\dim(SU_{12}/Sp_6) = 65$, this gives $\codim(\Sigma) \geq 31 > 20 = 4r$.

If $r = 6$, then $\dim(\Sigma) \leq \dim(E_6/F_4) + \dim(SU_4) = 41$. Since $\dim(SU_{14}/Sp_7) = 90$, this gives $\codim(\Sigma) \geq 49 > 24 = 4r$.

If $r \geq 7$, then $\dim(\Sigma) \leq \dim(E_6/F_4) + \dim(SU_{2(r-3)}/Sp_{r-3}) = 26 + (r-4)(2r-5)$. Since $\dim(SU_{2r+2}/Sp_{r+1}) = r(2r+3)$, this gives $\codim(\Sigma) \geq r(2r+3) - 26 - (r-4)(2r-5) = 14r - 46 > 4r$.

{\sc Case 3.} $\tilde{\Sigma}$ is a product of factors of the form $SU_{k_i+1}/SO_{k_i+1}$, $SU_{k_i+1}$ or $SU_{2k_i+2}/Sp_{k_i+1}$ ($2 \leq k_i \leq r$) with $k_1 + \ldots + k_b \leq r$. 

First assume that $\rk(\Sigma) = \rk(M) = r$. Let $\ca \subset T_o\Sigma \subset T_oM$ be a maximal abelian subspace. Then, by Proposition \ref{cardinal} and its proof, the root system $\tilde{\Phi}$ associated with the maximal abelian subspace $\ca$ of $T_o\Sigma$ is a root subsystem of the root system $\Phi$ associated with the maximal abelian subspace $\ca$ of $T_oM$. From (\ref {3}) it follows that there exists $0 \neq z \in \ca$ such that the Weyl group $\tilde{W}$ associated with $\tilde{\Phi}$ fixes $z$. This implies that $\Sigma$ is nonsemisimple, which contradicts the assumption that $\Sigma$ is semisimple.

Thus we have $\rk(\Sigma) < \rk(M) = r$. From the particular product form of $\tilde{\Sigma}$ it follows easily that $\dim(\Sigma) = \dim(\tilde{\Sigma}) \leq \dim(SU_{2r}/Sp_r) = (r-1)(2r+1)$ and therefore $\codim(\Sigma) = \dim(M) - \dim(\Sigma) \geq r(2r+3) - (r-1)(2r+1)  = 4r + 1 > 4r$.
\end{proof}

We can now state the main result of this section. 

\begin{thm} \label{LagGras}
For $M = SU_{2r+2}/Sp_{r+1}$ we have $i(M) = 4r = i_r(M)$ if $r \geq 4$ and $i(M) = 11 = i_r(M)$ if $r = 3$.
\end{thm}

\begin{proof}
We already know that $i_r(M) = 4r$ if $r \geq 4$ and $i_r(M) = 11$ if $r = 3$. Let $\Sigma$ be a maximal totally geodesic submanifold of $M$. If $\Sigma$ is nonsemisimple, then we have $\codim(\Sigma) \geq 4r$ by \cite[Theorem 4.2]{BO16}. Assume that $\Sigma$ is semisimple. If $\Sigma$ is locally reducible, it follows from Corollary \ref{corSigmaA_r} that every local de Rham factor of $\Sigma$ must have a root system of type $A$, which then implies $\codim(\Sigma) \geq 4r$ by Proposition \ref{factorsAAA}. Thus we can assume that $\Sigma$ is locally irreducible. If the root system of $\Sigma$ is not of type $A$, then we have three possibilities by Lemma \ref{AAA}:

(i): $\Sigma = G'/K'$ is inner, that is, $\rk(G') = \rk(K')$. Since $M = SU_{2r+2}/Sp_{r+1}$ is an outer symmetric space, then $\Sigma$ is a reflective submanifold by Proposition \ref{innSigma} and therefore $\codim(\Sigma) \geq i_r(M)$. 

(ii): $\max\{m_1, m_2\}\leq 2$. Then $\codim(\Sigma) \geq 4r$ by Proposition \ref{mleq2}.

(iii): $\Sigma = SO_{2s+n}/SO_sSO_{s+n}$ with $s \geq 3$ odd and $n \geq 4$ even. Then $\codim(\Sigma) \geq 4r$ by Proposition \ref{grassmannian}.

If the root system of $\Sigma$ is of type $A$, then $\codim(\Sigma) > 4r$ by Proposition \ref{factorsAAA}.
\end{proof}

\section{Quaternionic Grassmannians}

In this section we determine the index of the quaternionic Grassmann manifold $M = Sp_{2r+k}/Sp_rSp_{r+k}$, $r \geq 1$, $k \geq 0$. We already know the index for some values of $r$ and $k$ from results in \cite{BO16} (Table 4, Corollary 7.2 and Corollary 7.7). More precisely:

For $r = 1$ and $k = 0$ we have $i(M) = 1 (\neq 4r)$.

For $r = 1$ and $k > 0$ we have $i(M) = 4 (= 4r)$.

For $r = 2$ and $k = 0$ we have $i(M) = 6 (\neq 4r)$.

For $r \geq 2$ and $k \geq r-1$ we have $i(M) = 4r$.

It therefore remains to determine the index for $Sp_{2r+k}/Sp_rSp_{r+k}$ with $r \geq 3$ and $0 \leq k \leq r-2$.

\begin{lm} \label{spreduction}
Let $r \geq 3$. If $i(Sp_{2r}/Sp_rSp_r) = 4r$, then $i(Sp_{2r+k}/Sp_rSp_{r+k}) = 4r$ for all $k \geq 0$.
\end{lm}

\begin{proof}
The canonical inclusion $Sp_{2r} \subset Sp_{2r+k}$ leads to a canonical totally geodesic embedding of $\Sigma = Sp_{2r}/Sp_rSp_r = Sp_{2r} \cdot o$ into $M = Sp_{2r+k}/Sp_rSp_{r+k}$. From Proposition \ref{indexSigmaM} we obtain $4r = i(\Sigma) \leq i(M)$. On the other hand, from \cite{BO16} we know that $i_r(M) = 4r$. Since we always have $i(M) \leq i_r(M)$, we obtain $4r = i(\Sigma) \leq i(M) \leq i_r(M) = 4r$ and hence $i(M) = 4r$.
\end{proof}

It follows from Lemma \ref{spreduction} that it suffices to prove $i(Sp_{2r}/Sp_rSp_r) = 4r$ for $r \geq 3$. As we mentioned above, this equality does not hold for $r \in \{1,2\}$. We will prove this equality first for $r \in \{3,4,5\}$ and  then for arbitrary $r \geq 6$ by an inductive argument.

The next result provides useful bounds for the index of $Sp_{2r}/Sp_rSp_r$.

\begin{lm} \label{lowbound}
We have $4r-4 \leq i(Sp_{2r}/Sp_rSp_r) \leq 4r$ for all $r \geq 3$.
\end{lm}

\begin{proof}
The second inequality follows from the fact that the reflective index of $M = Sp_{2r}/Sp_rSp_r$ is equal to $4r$ for $r \geq 3$ (see \cite{BO16}). For the first inequality, consider the action of the isotropy group $K = Sp_rSp_r$ on the quaternionic Grassmannian $Sp_{2r}/Sp_rSp_r \cong G_r(\HH^{2r})$ of $r$-dimensional quaternionic subspaces of $\HH^{2r}$. This action induces a decomposition $\HH^{2r} = \HH^r \times \HH^r$. Define an $r$-dimensional quaternionic subspace $V$ of $\HH^{2r} = \HH^r \times \HH^r$ by $V = \{ (z,z) : z \in \HH^r \}$. The isotropy group of $K$ at $V$ is the diagonal subgroup $\Delta Sp_r$ and therefore the orbit $\Sigma = K \cdot V$ of $K$ containing $V$ is isometric to $Sp_r \cong Sp_rSp_r/\Delta Sp_r$. One can show that $\Sigma$ is a totally geodesic submanifold of $M$. In fact, if we consider the base point $o \in M$ as an $r$-dimensional quaternionic subspace of $\HH^{2r}$, then the orthogonal complement $o^\perp$ of $o$ in $\HH^{2r}$ is also a fixed point of the $K$-action on $M$. Thus $o^\perp$ is a pole of $o$ in $M$. The orbit of $K$ through the midpoint of a geodesic in $M$ connecting $o$ and $o^\perp$ is $\Sigma$. Thus $\Sigma \cong Sp_r$ is a centrosome of $M$ and therefore totally geodesic in $M$ (see \cite{NT95} for details). In \cite{BO17} we proved that $i(Sp_r) = 4r-4$. Using Proposition \ref{indexSigmaM} we then obtain $4r-4 = i(Sp_r) = i(\Sigma) \leq i(M)$.
\end{proof}

\begin{lm} \label{isotropybounds}
Let $\Sigma = G'/K'$ be a totally geodesic submanifold of $M = Sp_{2r}/Sp_rSp_r$ and $H$ be a (locally) irreducible factor of $K'$. Then the following inequalities hold:
\begin{align}
\rk(H) & \leq r, \label{bound1} \\
\dim(H) & \leq 2r^2 - 3r + 4 \label{bound2} \mbox{ if $\ch \not\cong \cs\cp_r$},\\
\rk(K') & \leq 2r, \label{bound3} \\
\dim(K') & < 2r(2r+1). \label{bound4}
\end{align}
\end{lm}

\begin{proof}
We prove these inequalities on Lie algebra level. Denote by $\pi$ the projection from $\ch$ into one of the two $\cs\cp_r$-factors. Then $\pi(\ch)$ is isomorphic to $\ch/\ker(\pi)$. Since $\ch$ is simple and $\pi(\ch) \neq \{0\}$, then $\pi(\ch)$ is isomorphic to $\ch$ and $\rk(\ch) = \rk(\pi(\ch)) \leq \rk(\cs\cp_r) = r$. Since $\pi(\ch) \neq \{0\}$ for at least one of the two projections, we proved (\ref{bound1}). If $\ch \not\cong \cs\cp_r$, then $\dim(\ch) = \dim(H) \leq \dim(Sp_r) - i_r(Sp_r) = (2r^2+r) - 4(r-1) = 2r^2-3r+4$, since the subgroup of $Sp_r$ with Lie algebra $\pi(\ch)$ is a totally geodesic submanifold of $Sp_r$ and $i_r(Sp_r) = 4(r-1)$ by \cite{BO17}. Since $\ck'$ is a subalgebra of $\cs\cp_r \oplus \cs\cp_r$, we obviously have $\rk(\ck') \leq \rk(\cs\cp_r \oplus \cs\cp_r) = 2r$ and $\dim(\ck') < \dim(\cs\cp_r \oplus \cs\cp_r) = 2(2r^2 + r) = 2r(2r+1)$.
\end{proof}

\begin{prop} \label{isp6sp3sp3}
For $M = Sp_6/Sp_3Sp_3$ we have $i(M) = 12$.
\end{prop}

\begin{proof}
We know from Table 4 in \cite{BO16} that $i_r(M) = 12$ and that $Sp_5/Sp_2Sp_3$ is a reflective submanifold of $M$ whose codimension is equal to $12$. Assume that there exists a maximal totally geodesic submanifold $\Sigma$ of $M$ with $d = \codim(\Sigma) < 12$. From Lemma \ref{lowbound} we obtain $d \in \{8,9,10,11\}$. We can slightly improve this. From Theorem 4.2 in \cite{BO16} we know that $\Sigma$ must be semisimple. The classification in \cite{BO16} of symmetric spaces with index $\leq 6$ tells us that $d \geq 7$. From Proposition 7.4 in \cite{BO16} we then get that $d(d-1) \geq 2(\dim(M) - \rk(M) - 1) = 64$, that is, $d \geq 9$. It follows that $d \in \{9,10,11\}$, or equivalently, $\dim(\Sigma) \in \{25,26,27\}$.

We write $\Sigma = G'/K'$ with $\ck' = [T_o\Sigma,T_o\Sigma]$ and $\cg' = [T_o\Sigma,T_o\Sigma] \oplus T_o\Sigma$.
Let $H$ be a (locally) irreducible factor of $K'$. From Lemma \ref{isotropybounds} we know that
\begin{align}
\rk(\ch) & \leq 3 \label{3bound1} ,\\
\rk(\ck') & \leq 6 \label{3bound3} , \\
\dim(\ck') & < 42 \label{3bound4} .
\end{align}

{\sc Case 1: $\rk(\Sigma) = 1$.} Then $\tilde\Sigma \in \{S^{25},S^{26},S^{27},\CC P^{13}\}$ and so $\ck' \in \{\cs\co_{25},\cs\co_{26},\cs\co_{27},\cu_{13}\}$. In all cases we have $\rk(\ck') > 6$, which contradicts (\ref{3bound3}).

{\sc Case 2: $\rk(\Sigma) = 2$.} By Theorem \ref{rankonefactors}, $\tilde\Sigma$ is irreducible or $\tilde\Sigma = S^{k_1} \times S^{k_2}$ with $k_1 \geq k_2 \geq 2$. If $\tilde\Sigma = S^{k_1} \times S^{k_2}$, then $\dim(\ck') = \dim(\cs\co_{k_1}) + \dim(\cs\co_{k_2}) = \frac{1}{2}(k_1(k_1-1) + k_2(k_2-1))$. Since $\dim(\Sigma) \in \{25,26,27\}$, we have $k_1 \geq 13$ and thus $\dim(\ck') \geq 78$, which contradicts (\ref{3bound4}). Thus $\tilde\Sigma$ is irreducible.
Since $\dim(\Sigma) \in \{25,26,27\}$ and $\rk(\Sigma) = 2$, we have only two possibilities, namely $\tilde\Sigma = SO_{15}/SO_2SO_{13}$ and $\tilde\Sigma = E_6/F_4$. Since $\rk(\cs\co_{13}) = 6$ and $\rk(\cf_4) = 4$, we can exclude both possibilities using (\ref{3bound1}).

{\sc Case 3: $\rk(\Sigma) = 3$.} By Theorem \ref{rankonefactors}, $\tilde\Sigma$ is irreducible. Since $\dim(\Sigma) \in \{25,26,27\}$ and $\rk(\Sigma) = 3$, we have only one possibility, namely $\tilde\Sigma = SO_{12}/SO_3SO_9$. Since $\rk(\cs\co_9) = 4$, we can exclude this possibility using (\ref{3bound1}).

Altogether it now follows that there exists no maximal totally geodesic submanifold $\Sigma$ of $M$ with $\codim(\Sigma) < 12 = i_r(M)$, and therefore $i(M) = 12$.
\end{proof}

\begin{prop} \label{isp8sp4sp4}
For $M = Sp_8/Sp_4Sp_4$ we have $i(M) = 16$.
\end{prop}

\begin{proof}
We know from Table 4 in \cite{BO16} that $i_r(M) = 16$ and that $Sp_7/Sp_3Sp_4$ is a reflective submanifold of $M$ whose codimension is equal to $16$. Assume that there exists a maximal totally geodesic submanifold $\Sigma$ of $M$ with $d = \codim(\Sigma) < 16$. From Lemma \ref{lowbound} we obtain $d \in \{12,13,14,15\}$, or equivalently, $\dim(\Sigma) \in \{49,50,51,52\}$. 

We write $\Sigma = G'/K'$ with $\ck' = [T_o\Sigma,T_o\Sigma]$ and $\cg' = [T_o\Sigma,T_o\Sigma] \oplus T_o\Sigma$.
Let $H$ be a (locally) irreducible factor of $K'$. From Lemma \ref{isotropybounds} we know that
\begin{align}
\rk(\ch) & \leq 4 \label{4bound1} ,\\
\dim(\ch) & \leq 24 \label{4bound2} \mbox{ if $\ch \not\cong \cs\cp_4$} ,\\
\rk(\ck') & \leq 8 \label{4bound3} , \\
\dim(\ck') & < 72 \label{4bound4} .
\end{align}

{\sc Case 1: $\rk(\Sigma) = 1$.} Then $\tilde\Sigma \in \{S^{49},S^{50},S^{51},S^{52},\CC P^{25},\CC P^{26},\HH P^{13}\}$ and so $\ck' \in \{\cs\co_{49},\cs\co_{50},\cs\co_{51},\cs\co_{52},\cu_{25},\cu_{26},\cs\cp_{13}\oplus
\cs\cp_1\}$. In all cases we have $\rk(\ck') > 8$, which contradicts (\ref{4bound3}).

{\sc Case 2: $\rk(\Sigma) = 2$.} By Theorem \ref{rankonefactors}, $\tilde\Sigma$ is irreducible or $\tilde\Sigma = S^{k_1} \times S^{k_2}$ with $k_1 \geq k_2 \geq 2$. If $\tilde\Sigma = S^{k_1} \times S^{k_2}$, then $\dim(\ck') = \dim(\cs\co_{k_1}) + \dim(\cs\co_{k_2}) = \frac{1}{2}(k_1(k_1-1) + k_2(k_2-1))$. Since $\dim(\Sigma) \in \{49,50,51,52\}$, we have $k_1 \geq 25$ and thus $\dim(\ck') \geq 300$, which contradicts (\ref{4bound4}). Thus $\tilde\Sigma$ is irreducible. Since $\rk(\Sigma) = 2$ and $\dim(\Sigma) \in \{49,50,51,52\}$, we have only three possibilities, namely $\tilde\Sigma = SO_{27}/SO_2SO_{25}$, $\tilde\Sigma = SO_{28}/SO_2SO_{26}$ and $\tilde\Sigma = SU_{15}/S(U_2U_{13})$. In all cases we have $\rk(\ck') > 8$, which contradicts (\ref{4bound3}).

{\sc Case 3: $\rk(\Sigma) = 3$.} By Theorem \ref{rankonefactors}, $\tilde\Sigma$ is irreducible. Since $\rk(\Sigma) = 3$ and $\dim(\Sigma) \in \{49,50,51,52\}$, we have only one possibility, namely $\tilde\Sigma = SO_{20}/SO_3SO_{17}$. In this case we have $\rk(\ck') > 8$, which contradicts (\ref{4bound3}).

{\sc Case 4: $\rk(\Sigma) = 4$.} By Theorem \ref{rankonefactors}, $\tilde\Sigma$ is irreducible or the product of two symmetric spaces of rank $2$. Firstly, assume that $\tilde\Sigma$ is irreducible. Since $\dim(\Sigma) \in \{49,50,51,52\}$ and $\rk(\Sigma) = 4$, we have only one possibility, namely $\tilde\Sigma = SO_{17}/SO_4SO_{13}$. In this case we have $\ch = \cs\co_{13}$ and $\rk(\ch) = 6$, which contradicts (\ref{4bound1}). Next, assume that $\tilde\Sigma = \tilde\Sigma_1 \times \tilde\Sigma_2$, where $\tilde\Sigma_1$ and $\tilde\Sigma_2$ are irreducible symmetric spaces of rank $2$. we can assume that $\dim(\tilde\Sigma_1) \geq \dim(\tilde\Sigma_2)$. Then $25 \leq \dim(\tilde\Sigma_1) \leq 47$, using the fact that $5$ is the lowest dimension of an irreducible symmetric space of rank $2$. We discuss the various possibilities.

If $\tilde\Sigma_1 = SO_{2+k}/SO_2SO_k$, $13 \leq k \leq 23$, then $\ch = \cs\co_k$ with $k \in \{13,\ldots,23\}$ and thus $\rk(\ch) > 4$, which contradicts (\ref{4bound1}).

If $\tilde\Sigma_1 = SU_{2+k}/S(U_2U_k)$, $7 \leq k \leq 11$, then $\ch = \cs\cu_k$ with $k \in \{7,\ldots,11\}$ and thus $\rk(\ch) > 4$, which contradicts (\ref{4bound1}).

If $\tilde\Sigma_1 = Sp_{2+k}/Sp_2Sp_k$, $4 \leq k \leq 5$. For $k = 5$ we have $\ch = \cs\cp_5$ and thus $\rk(\ch) > 4$, which contradicts (\ref{4bound1}). If $k = 4$, then $\dim(\tilde\Sigma_1) = 32$ and thus $17 \leq \dim(\tilde\Sigma_2) \leq 20$. Moreover, (\ref{4bound3}) implies that the rank of the isotropy group of $\tilde\Sigma_2$ must be $\leq 2$. It is easy to check that there does not exist an irreducible symmetric space of rank $2$ with these properties. 

If $\tilde\Sigma_1 = E_6/F_4$, then $\ch = \cf_4$ and $\dim(\ch) = 52$, which contradicts (\ref{4bound2}).

If $\tilde\Sigma_1 = E_6/Spin_{10}U_1$, then $\ch = \cs\co_{10}$ and $\rk(\ch) = 5$, which contradicts (\ref{4bound1}).

Altogether it now follows that there exists no maximal totally geodesic submanifold $\Sigma$ of $M$ with $\codim(\Sigma) < 16 = i_r(M)$, and therefore $i(M) = 16$.
\end{proof}

\begin{prop} \label{isp10sp5sp5}
For $M = Sp_{10}/Sp_5Sp_5$ we have $i(M) = 20$.
\end{prop}

\begin{proof}
We know from Table 4 in \cite{BO16} that $i_r(M) = 20$ and that $Sp_9/Sp_4Sp_5$ is a reflective submanifold of $M$ whose codimension is equal to $20$. Assume that there exists a maximal totally geodesic submanifold $\Sigma$ of $M$ with $d = \codim(\Sigma) < 20$. From Lemma \ref{lowbound} we obtain $d \in \{16,17,18,19\}$, or equivalently, $\dim(\Sigma) \in \{81,82,83,84\}$. 

We write $\Sigma = G'/K'$ with $\ck' = [T_o\Sigma,T_o\Sigma]$ and $\cg' = [T_o\Sigma,T_o\Sigma] \oplus T_o\Sigma$.
Let $H$ be a (locally) irreducible factor of $K'$. From Lemma \ref{isotropybounds} we know that
\begin{align}
\rk(\ch) & \leq 5 \label{5bound1} ,\\
\dim(\ch) & \leq 39 \label{5bound2} \mbox{ if $\ch \not\cong \cs\cp_5$},\\
\rk(\ck') & \leq 10 \label{5bound3} , \\
\dim(\ck') & < 110 \label{5bound4} .
\end{align}

{\sc Case 1: $\rk(\Sigma) = 1$.} Then $\tilde\Sigma \in \{S^{81},S^{82},S^{83},S^{84},\CC P^{41},\CC P^{42},\HH P^{21}\}$ and so $\ck' \in \{\cs\co_{81},\cs\co_{82},\cs\co_{83},\cs\co_{84},\cu_{41},\cu_{42},\cs\cp_{21}\oplus
\cs\cp_1\}$. In all cases we have $\rk(\ck') > 10$, which contradicts (\ref{5bound3}).

{\sc Case 2: $\rk(\Sigma) = 2$.} By Theorem \ref{rankonefactors}, $\tilde\Sigma$ is irreducible or $\tilde\Sigma = S^{k_1} \times S^{k_2}$ with $k_1 \geq k_2 \geq 2$. If $\tilde\Sigma = S^{k_1} \times S^{k_2}$, then $\dim(\ck') = \dim(\cs\co_{k_1}) + \dim(\cs\co_{k_2}) = \frac{1}{2}(k_1(k_1-1) + k_2(k_2-1))$. Since $\dim(\Sigma) \in \{81,82,83,84\}$, we have $k_1 \geq 41$ and thus $\dim(\ck') \geq 820$, which contradicts (\ref{5bound4}). Thus $\tilde\Sigma$ is irreducible. Since $\rk(\Sigma) = 2$ and $\dim(\Sigma) \in \{81,82,83,84\}$, we have only three possibilities, namely $\tilde\Sigma = SO_{43}/SO_2SO_{41}$, $\tilde\Sigma = SO_{44}/SO_2SO_{42}$ and $\tilde\Sigma = SU_{23}/S(U_2U_{21})$. In all cases we have $\rk(\ck') > 10$, which contradicts (\ref{5bound3}).

{\sc Case 3: $\rk(\Sigma) = 3$.} By Theorem \ref{rankonefactors}, $\tilde\Sigma$ is irreducible. Since $\rk(\Sigma) = 3$ and $\dim(\Sigma) \in \{81,82,83,84\}$, we have four possibilities, namely $\tilde\Sigma = SO_{30}/SO_3SO_{27}$, $\tilde\Sigma = SO_{31}/SO_3SO_{28}$, $\tilde\Sigma = SU_{17}/S(U_3U_{14})$ and $\tilde\Sigma = Sp_{10}/Sp_3Sp_7$. In the first three cases we have $\rk(\ck') > 10$, which contradicts (\ref{5bound3}). In the last case we can choose $\ch = \cs\cp_7$, then $\rk(\ch) = 7$, which contradicts (\ref{5bound1}).

{\sc Case 4: $\rk(\Sigma) = 4$.} By Theorem \ref{rankonefactors}, $\tilde\Sigma$ is irreducible or the product of two symmetric spaces of rank $2$. Firstly, assume that $\tilde\Sigma$ is irreducible. Since $\dim(\Sigma) \in \{81,82,83,84\}$ and $\rk(\Sigma) = 4$, we have only one possibility, namely $\tilde\Sigma = SO_{25}/SO_4SO_{21}$. In this case we can choose $\ch = \cs\co_{21}$ and so $\rk(\ch) = 10$, which contradicts (\ref{5bound1}). Next, assume that $\tilde\Sigma = \tilde\Sigma_1 \times \tilde\Sigma_2$, where $\tilde\Sigma_1$ and $\tilde\Sigma_2$ are irreducible symmetric spaces of rank $2$. we can assume that $\dim(\tilde\Sigma_1) \geq \dim(\tilde\Sigma_2)$. Then $41 \leq \dim(\tilde\Sigma_1) \leq 79$, using the fact that $5$ is the lowest dimension of an irreducible symmetric space of rank $2$. We discuss the various possibilities.

If $\tilde\Sigma_1 = SO_{2+k}/SO_2SO_k$, $21 \leq k \leq 39$, then $\ch = \cs\co_k$ with $k \in \{21,\ldots,39\}$ and thus $\rk(\ch) > 5$, which contradicts (\ref{5bound1}).

If $\tilde\Sigma_1 = SU_{2+k}/S(U_2U_k)$, $11 \leq k \leq 19$, then $\ch = \cs\cu_k$ with $k \in \{11,\ldots,19\}$ and thus $\rk(\ch) > 5$, which contradicts (\ref{5bound1}).

If $\tilde\Sigma_1 = Sp_{2+k}/Sp_2Sp_k$, $6 \leq k \leq 9$, then $\ch = \cs\cp_k$ with $k \in \{6,\ldots,9\}$ and thus $\rk(\ch) > 5$, which contradicts (\ref{5bound1}). 

{\sc Case 5: $\rk(\Sigma) = 5$.} By Theorem \ref{rankonefactors}, $\tilde\Sigma$ is irreducible or the product of a symmetric space of rank $2$ and a symmetric space of rank $3$. However, there are no irreducible symmetric paces of rank $5$ and dimension in $\{81,82,83,84\}$. Consequently, $\tilde\Sigma = \tilde\Sigma_1 \times \tilde\Sigma_2$, where $\tilde\Sigma_1$ and $\tilde\Sigma_2$ are irreducible symmetric spaces of rank 2 or 3 and $\rk(\tilde\Sigma_1) + \rk(\tilde\Sigma_2) = 5$. We can assume that $\dim(\tilde\Sigma_1) \geq \dim(\tilde\Sigma_2)$. Then $41 \leq \dim(\tilde\Sigma_1) \leq 79$, using the fact that $5$ is the lowest dimension of an irreducible symmetric space of rank $\geq 2$. If $\rk(\tilde\Sigma_1) = 2$, we can use the arguments given in the previous Case 4. Assume that $\rk(\tilde\Sigma_1) = 3$. We discuss the various possibilities.

If $\tilde\Sigma_1 = SO_{3+k}/SO_3SO_k$, $14 \leq k \leq 26$, then $\ch = \cs\co_k$ with $k \in \{14,\ldots,26\}$ and thus $\rk(\ch) > 5$, which contradicts (\ref{5bound1}).

If $\tilde\Sigma_1 = SU_{3+k}/S(U_3U_k)$, $7 \leq k \leq 13$, then $\ch = \cs\cu_k$ with $k \in \{7,\ldots,13\}$ and thus $\rk(\ch) > 5$, which contradicts (\ref{5bound1}).

If $\tilde\Sigma_1 = Sp_{3+k}/Sp_3Sp_k$, $4 \leq k \leq 6$, we need different arguments.

If $k = 6$, then we can choose $\ch = \cs\cp_6$ and thus $\rk(\ch) > 5$, which contradicts (\ref{5bound1}). 

If $k = 5$, then $\dim(\tilde\Sigma_1) = 60$ and thus $21 \leq \dim(\tilde\Sigma_2) \leq 24$. It follows that $\tilde\Sigma_2 = SO_{13}/SO_2SO_{11}$, or $\tilde\Sigma_2 = SO_{14}/SO_2SO_{12}$, or $\tilde\Sigma_2 = SU_8/S(U_2U_6)$, or $\tilde\Sigma_2 = Sp_5/Sp_2Sp_3$. Since the isotropy algebra of $\tilde\Sigma_1$ has rank $8$, the isotropy algebra of $\tilde\Sigma_2$ must have rank $\leq 2$ by (\ref{5bound3}). However, in all four cases the isotropy algebra of $\tilde\Sigma_2$ has rank $> 2$, which gives a contradiction.

If $k = 4$, then $\dim(\tilde\Sigma_1) = 48$ and thus $33 \leq \dim(\tilde\Sigma_2) \leq 36$. It follows that $\tilde\Sigma_2 = SO_{19}/SO_2SO_{17}$, or $\tilde\Sigma_2 = SO_{20}/SO_2SO_{18}$, or $\tilde\Sigma_2 = SU_{11}/S(U_2U_9)$. Since the isotropy algebra of $\tilde\Sigma_1$ has rank $7$, the isotropy algebra of $\tilde\Sigma_2$ must have rank $\leq 3$ by (\ref{5bound3}). However, in all three cases the isotropy algebra of $\tilde\Sigma_2$ has rank $> 3$, which gives a contradiction.

If $\tilde\Sigma_1 = SO_{14}/U_7$, then $\ch = \cs\cu_7$ satisfies $\rk(\ch) = 6$, which contradicts (\ref{5bound1}).

If $\tilde\Sigma_1 = E_7/E_6U_1$, then $\ch = \ce_6$ satisfies $\rk(\ch) = 6$, which contradicts (\ref{5bound1}).

Altogether it now follows that there exists no maximal totally geodesic submanifold $\Sigma$ of $M$ with $\codim(\Sigma) < 20 = i_r(M)$, and therefore $i(M) = 20$.
\end{proof}

\begin{prop} \label{isp2rsprspr}
For $M = Sp_{2r}/Sp_rSp_r$ ($r \geq 3$) we have $i(M) = 4r$.
\end{prop}

\begin{proof}
Consider the symmetric space $N = Sp_{2(r+3)}/Sp_{r+3}Sp_{r+3}$. We already know from \cite{BO16} that $i_r(N) = 4(r+3)$. The symmetric space $N^+ = Sp_{r+3}/Sp_3Sp_r \times Sp_{r+3}/Sp_rSp_3$ is a polar of $N$ with corresponding meridian $N^- = Sp_{2r}/Sp_rSp_r \times Sp_6/Sp_3Sp_3$ (see \cite{CN78} for details). Both $N^+$ and $N^-$ are reflective submanifolds of $N$ and $o \in N^-$.

Let $\Sigma$ be a totally geodesic submanifold of $N$ with $\codim(\Sigma) = i(N)$ and $o \in \Sigma$. Then $\Sigma' = \Sigma \cap N^-$ is a totally geodesic submanifold of $N^-$. By construction, the codimension $\codim_N(\Sigma)$ of $\Sigma$ in $N$ and the codimension $\codim_{N^-}(\Sigma')$ of $\Sigma'$ in $N^-$ satisfy
\[
\codim_N(\Sigma) \geq \codim_{N^-}(\Sigma').
\]
We define subspaces $V_1,V_2,W$ of $T_oN$ by $V_1 = T_o(Sp_{2r}/Sp_rSp_r)$, $V_2 = T_o(Sp_6/Sp_3Sp_3)$ and $W = T_o\Sigma$. Then, by construction, $T_oN^- = V_1 \oplus V_2$.

For $j \in \{1,2\}$ we define $K^j = \{k \in K : V_j \subseteq d_ok(W)\}$, where $K = Sp_{r+3}Sp_{r+3}$. If $k \in K \setminus K^j$, then $V_j$ is not a subspace of $d_ok(W)$. By continuity of the action of $K$ on $T_oN$, there exists an open neighborhood ${\mathcal U}_j$ of $k$ in $K$ such that $V_j$ is not a subspace of $d_ok(W)$ for all $k \in {\mathcal U}_j$. Thus $K \setminus K^j$ is an open subset of $K$.

Let $0 \neq u \in T_oN$, $0 \neq v_j \in V_j$, and assume that $u$ is perpendicular to $W$. We define the analytic function 
\[
f_{v_j,u} : K \to \RR\ , k \mapsto \langle v_j, d_ok(u) \rangle. 
\] 
Assume that $K^j$ contains a nonempty open subset $\Omega$ of $K$. Then $f_{v_j,u}|_\Omega = 0$ and thus $f_{v_j,u} = 0$ by analyticity of $f_{v_j,u}$ and since $K$ is connected. It follows that $\RR v_j \subseteq d_ok(W)$ for all $k \in K$. This is a contradiction since $\bigcap_{k \in K} d_ok(W) = \{0\}$. It follows that $K^j$ does not contain any nonempty open subsets of $K$. 

Altogether we now see that $K \setminus K^j = \{ k \in K : V_j \not\subseteq d_ok(W)\}$ is an open and dense subset of $K$. It follows that 
\[
(K \setminus K^1) \cap (K \setminus K^2) = \{ k \in K : V_1 \not\subseteq d_ok(W) \mbox{ and } V_2 \not\subseteq d_ok(W) \}
\]
is an open and dense subset of $K$. We can therefore assume, without loss of generality, that $V_1$ and $V_2$ are not contained in $T_o\Sigma'$. Using Proposition \ref{productestimate} we then obtain
\begin{align*}
i(Sp_{2(r+3)}/Sp_{r+3}Sp_{r+3}) & = \codim_N(\Sigma) \geq \codim_{N^-}(\Sigma') \\
& \geq  i(Sp_{2r}/ Sp_rSp_r) + i(Sp_6/ Sp_3Sp_3).
\end{align*}
By induction over $r$ we can now prove the assertion. For $r \in \{3,4,5\}$ we already know that $i(Sp_{2r}/ Sp_rSp_r)  = 4r$ by Propositions \ref{isp6sp3sp3}, \ref{isp8sp4sp4} and \ref{isp10sp5sp5}. For $r \geq 6$ we can then use the previous inequality and the induction hypothesis and obtain
\begin{align*}
4(r+3) & = i_r(Sp_{2(r+3)}/Sp_{r+3}Sp_{r+3}) \\
& \geq i(Sp_{2(r+3)}/Sp_{r+3}Sp_{r+3})   \\
& \geq  i(Sp_{2r}/ Sp_rSp_r) + i(Sp_6/ Sp_3Sp_3) = 4r + 12 = 4(r+3).
\end{align*}
This finishes the proof.
\end{proof}

\begin{re}
\rm Proposition \ref{productestimate} can be generalized to the case that $M = M_1 \times \ldots \times M_k$ is the Riemannian product of a finite number of irreducible factors. Using similar arguments as in the proof of Proposition \ref{isp2rsprspr}, we can then show the following: 

Let $\Sigma = \Sigma_1\times \ldots \times \Sigma_k$ be a totally geodesic submanifold of a Riemannian symmetric space $M$ of noncompact type, where $\Sigma_1, \ldots, \Sigma _k$ are irreducible factors of $\Sigma$. Then 
\[
i(\Sigma _1) + \ldots + i(\Sigma _k) \leq i(M).
\] 
In particular, if $i(M) = \rk(M)$, then $i(\Sigma_1) + \ldots + i(\Sigma _k) \leq \rk(M)$. Therefore, if $\rk(\Sigma) = \rk(M)$, then $i(\Sigma _\nu) = \rk(\Sigma _\nu)$ for all $\nu \in \{1, \ldots , k\}$. The possible factors are then known from \cite{BO16}.
\end{re}

\begin{thm} \label{quatGras}
We have $i(Sp_{2r+k}/Sp_rSp_{r+k}) = 4r$ for all $r \geq 3$ and  $k \geq 0$.
\end{thm}

\begin{proof}
This now follows from Lemma \ref{spreduction} and Proposition \ref{isp2rsprspr}.
\end{proof}

\section{Irreducible Hermitian symmetric spaces}

In this section we study the index of irreducible Hermitian symmetric spaces. Our first result states that a maximal totally geodesic submanifold of sufficiently small codimension in an irreducible Hermitian symmetric space is a complex submanifold.

\begin{prop} \label{maxHerm}
Let $M = G/K$ be an irreducible Hermitian symmetric space. Every maximal totally geodesic submanifold $\Sigma$ of $M$ with $\codim(\Sigma) < \frac{1}{2}\dim(M)$ is a (semisimple) complex submanifold. 
\end{prop}

\begin{proof}
By duality, we can assume that $M$ is of noncompact type.

Let $\Sigma$ be a maximal totally geodesic submanifold of $M$ with $\codim(\Sigma) < \frac{1}{2}\dim(M)$. We can assume that $o \in \Sigma$. Let $\cg = \ck + \cp$ be the corresponding Cartan decomposition of $\cg$. The center $\cz(\ck)$ of $\ck$ is $1$-dimensional and there exists $Z \in \cz(\ck)$ such that $J = \ad(Z)$ is the complex structure on $\cp \cong T_oM$. The differential $d_oz$ at $o$ of the isometry $z = \Exp(\frac{\pi}{2}Z) \in K$ of $M$  is $d_oz : T_oM \to T_oM,\ X \mapsto JX$. Then $\Sigma^J = z(\Sigma)$ is also a maximal totally geodesic submanifold of $M$ with $o \in \Sigma^J$. By construction, we have $T_o\Sigma^J = JT_o\Sigma$. 

As usual, we write $\Sigma = G'/K'$ with $\cg' = \ck' + \cp' \subset \ck + \cp = \cg$, where $\cp' = T_o\Sigma$ and $\ck' = [\cp',\cp']$. Then, since $Z \in \cz(\ck)$, we have $\Sigma^J = G''/K'$ with $G'' = zG'z^{-1}$. Now consider the de Rham decomposition $\Sigma = \Sigma_0 \times \Sigma_1 \times \ldots \times \Sigma_l$ of $\Sigma$, where $\Sigma_0$ is the, possibly $0$-dimensional, Euclidean factor. Then, by construction, the de Rham decomposition of $\Sigma^J$ is $\Sigma^J = \Sigma_0^J \times \Sigma_1^J \times \ldots \times \Sigma_l^J$ with $\Sigma_i^J = z(\Sigma_i)$.

The intersection $\cq = T_o\Sigma \cap T_o\Sigma^J = \cp' \cap J\cp'$ is a $J$-invariant Lie triple system in $\cp$. Since $\dim(\Sigma^J) = \dim(\Sigma) > \frac{1}{2}\dim(M)$, we have $\dim(\cq) > 0$.  As both $\cp'$ and $J\cp'$ are $\Ad(K')$-invariant, the intersection $\cq =  \cp' \cap J\cp'$ is also $\Ad(K')$-invariant.

Since maximal flats of irreducible Hermitian symmetric spaces are totally real submanifolds, the Euclidean factor $\Sigma_0$ is a totally real submanifold of $M$ and therefore $T_o\Sigma_0 \cap T_o\Sigma_0^J = \{0\}$. Since $\Ad(K')$ acts irreducibly on each tangent space $T_o\Sigma_j$ (and $T_o\Sigma_j^J$) for $1 \leq j \leq l$, we see that $\cq = \bigoplus_{i \in I} T_o\Sigma_i$ for some nonempty subset $I$ of $\{1,\ldots,l\}$ and $T_o\Sigma_i^J = T_o\Sigma_i$ for all $i \in I$. 

Let $\cz_\cp(\cq) = \{ U \in \cp : [U,\cq] = \{0\}\}$ be the centralizer of $\cq$ in $\cp$ and put $\crr = \cz_\cp(\cq) + \cq$. We claim that $\crr$ is a Lie triple system in $\cp$ containing both $T_o\Sigma$ and $T_o\Sigma^J$. If $U,V,W \in \cz_\cp(\cq)$, then $[[U,V],\cq] = \{0\}$ by the Jacobi identity, which implies $[[[U,V],W],\cq] = \{0\}$ by the Jacobi identity. Thus $[[U,V],W] \in \cz_\cp(\cq)$, which shows that $\cz_\cp(\cq)$ is a Lie triple system in $\cp$. For $U,V,W \in \cz_\cp(\cq)$ and $U',V',W' \in \cq$ we get $[[U+U',V+V'],W+W'] = [[U,V],W]  + [[U',V'],W'] \in  \cz_\cp(\cq) + \cq = \crr$ by a straightforward calculation. Thus $\crr$ is a Lie triple system.  Since $\Sigma$ is a Riemannian product $\Sigma = \Sigma_0 \times \Sigma_1 \times \ldots \times \Sigma_l$, we have $[T_o\Sigma_i,T_o\Sigma_j] = \{0\}$ for all $0 \leq i < j \leq l$. As $\cq = \bigoplus_{i \in I} T_o\Sigma_i$, it follows that $T_o\Sigma_i \subseteq \cz_\cp(\cq)$ for all $i \notin I$. Similarly, we have $T_o\Sigma_i^J \subseteq \cz_\cp(\cq)$ for all $i \notin I$. Altogether we see that $\crr$ is a Lie triple system in $\cp$ containing $T_o\Sigma$ and $T_o\Sigma^J$. 

Assume that $\crr = \cp$. Then we can write $\crr = \cp = \cq \oplus \cq^\perp$ with $[\cq,\cq^\perp] = \{0\}$. This implies that $M$ is reducible, which is a contradiction. Thus $\crr$ is properly contained in $\cp$.

Since $\Sigma$ (and $\Sigma^J$) is a maximal totally geodesic submanifold of $M$, we must have $\crr = T_o\Sigma = T_o\Sigma^J$, which means that $\Sigma = \Sigma^J$ is a complex submanifold of $M$. In particular, the de Rham decomposition of $\Sigma$ has no Euclidean factor and thus $\Sigma$ is semisimple.
\end{proof}

Our next result states that a maximal totally geodesic submanifold with sufficiently small codimension in an irreducible Hermitian symmetric space must be a reflective submanifold.

\begin{prop} \label{estimate} 
Let $M$ be an irreducible Hermitian symmetric space of noncompact type with $\rk(M) \geq 2$ and let $\Sigma$ be a maximal totally geodesic submanifold of $M$ with $\codim(\Sigma) < \frac{1}{2}\dim(M)$. If 
\[
\textstyle{\frac{1}{4}}\codim(\Sigma)^2 + \codim(\Sigma) +\rk(\Sigma) < \dim(M),
\] 
then $\Sigma$ is a reflective submanifold of $M$. In particular, if 
\[
\textstyle{\frac{1}{4}}\codim(\Sigma)^2 + \codim(\Sigma)  < \dim(M) - \rk(M),
\] 
then $\Sigma$ is a reflective submanifold of $M$.
\end{prop}

\begin{proof}
By Proposition \ref{maxHerm}, $\Sigma$ is a semisimple complex submanifold of $M$. Therefore the slice representation $\rho : K' \to SO(\nu_o\Sigma)$ acts by unitary transformations.
As usual, we write $\Sigma = G'/K'$ with $\ck' = [T_o\Sigma,T_o\Sigma]$ and $\cg' = \ck' + T_o\Sigma$. Then $\rho(K') \subseteq U(\nu_o\Sigma)$. If $\dim(K') > \dim(U(\nu_o\Sigma)) = \frac{1}{4}\codim(\Sigma)^2$, then the kernel of the slice representation $\rho$ must have positive dimension and therefore $\Sigma$ is a reflective submanifold of $M$ by Proposition 3.4 in \cite{BO16}. A principal $K'$-orbit on $\Sigma$ has dimension $\dim(M)-\codim(\Sigma)-\rk(\Sigma)$ and thus $\dim(K') \geq \dim(M) - \codim(\Sigma) - \rk(\Sigma)$. Consequently, if $\frac{1}{4}\codim(\Sigma)^2 < \dim(M) - \codim(\Sigma) - \rk(\Sigma)$, then $\Sigma$ is a reflective submanifold of $M$. The last statement follows from the fact that $\rk(\Sigma) \leq \rk(M)$.
\end{proof}

We now apply the previous two results to irreducible Hermitian symmetric spaces $M$ with $r = \rk(M) \geq 2$. We put $d = \codim(\Sigma)$ and $n = \dim(M)$.

\smallskip
For $M = SO^o_{2,2+k}/SO_{2}SO_{2+k}$ ($k \geq 1$) we have $n = 2k+4$, $r = 2$ and $i_r(M) = 2$. By Proposition \ref{maxHerm}, a maximal totally geodesic submanifold $\Sigma$ of $M$ with $d < k+2$ must be complex. This immediately implies $i(M) = 2 = i_r(M)$.

\smallskip
For $M = Sp_r(\RR)/U_r$ we have $n = r^2+r$ and $i_r(M) = 2r-2$. By Proposition \ref{maxHerm}, a maximal totally geodesic submanifold $\Sigma$ of $M$ with $d < \frac{1}{2}(r^2+r)$ must be complex. Assume that $d  < i_r(M)$. Then $d  \in \{2,4,\ldots,2r-4\}$. We have 
\[
\frac{1}{4}d^2 + d \leq \frac{1}{4}(2r-4)^2 + (2r-4) =  r^2 -2r < r^2 = n-r.
\]
It follows from Proposition \ref{estimate} that $\Sigma$ is reflective, which contradicts $d < i_r(M)$. Thus we must have $i(M) = i_r(M)$.

\smallskip
For $M = SU_{r,r+k}/S(U_rU_{r+k})$ ($k \geq 0$) we have $n = 2r(r+k)$ and $i_r(M) = 2r$. By Proposition \ref{maxHerm}, a maximal totally geodesic submanifold $\Sigma$ of $M$ with $d < r(r+k)$ must be complex. Assume that $d  < i_r(M)$. Then $d  \in \{2,4,\ldots,2r-2\}$. We have 
\[
\frac{1}{4}d^2 + d \leq \frac{1}{4}(2r-2)^2 + (2r-2) =  r^2-1 < 2r^2+2rk-r = n-r. 
\]
It follows from Proposition \ref{estimate} that $\Sigma$ is reflective, which contradicts $d < i_r(M)$. Thus we must have $i(M) = i_r(M)$.

\smallskip
For $M = SO^*_{4r}/U_{2r}$ we have $n = 4r^2-2r$ and $i_r(M) = 4r-2$. By Proposition \ref{maxHerm}, a maximal totally geodesic submanifold $\Sigma$ of $M$ with $d < 2r^2-r$ must be complex. Assume that $d  < i_r(M)$. Then $d  \in \{2,4,\ldots,4r-4\}$. We have 
\[
\frac{1}{4}d^2 + d \leq \frac{1}{4}(4r-4)^2 + (4r-4) =  4r^2 -4r < 4r^2 - 3r = n-r. 
\]
It follows from Proposition \ref{estimate} that $\Sigma$ is reflective, which contradicts $d < i_r(M)$. Thus we must have $i(M) = i_r(M)$.

\smallskip
For $M = SO^*_{4r+2}/U_{2r+1}$ we have $n = 4r^2+2r$ and $i_r(M) = 4r$. By Proposition \ref{maxHerm}, a maximal totally geodesic submanifold $\Sigma$ of $M$ with $d < 2r^2+r$ must be complex. Assume that $d  < i_r(M)$. Then $d  \in \{2,4,\ldots,4r-2\}$. We have 
\[
\frac{1}{4}d^2 + d \leq \frac{1}{4}(4r-2)^2 + (4r-2) =  4r^2 - 1 < 4r^2 + r = n-r. 
\]
It follows from Proposition \ref{estimate} that $\Sigma$ is reflective, which contradicts $d < i_r(M)$. Thus we must have $i(M) = i_r(M)$.

\smallskip
From these calculations we conclude:

\begin{thm} \label{Hermsymm}
For an irreducible Hermitian symmetric space of classical type we have $i(M) = i_r(M)$.
\end{thm}

\begin{re} 
\rm It is somewhat surprising that this argument is inconclusive for the irreducible Hermitian symmetric spaces of exceptional type.

For $M = E_6^{-14}/Spin_{10}U_1$ we have $n = 32$, $r = 2$ and $i_r(M) = 12$. By Proposition \ref{maxHerm}, a maximal totally geodesic submanifold $\Sigma$ of $M$ with $d < 16$ must be complex. Assume that $d  < i_r(M)$. Then $d  \in \{2,4,6,8,10\}$. We have 
\[
\frac{1}{4}d^2 + d  <  30 = n-r \quad \Longleftrightarrow \quad d \in \{2,4,6,8\}. 
\]
It follows from Proposition \ref{estimate} that the index of $M$ is either $10$ or $12$. We proved in \cite{BOR19}, with different methods, that $i(M) = 12 = i_r(M)$.

For $M = E_7^{-25}/E_6U_1$ we have $n = 54$, $r = 3$ and $i_r(M) = 22$. By Proposition \ref{maxHerm}, a maximal totally geodesic submanifold $\Sigma$ of $M$ with $d < 27$ must be complex. Assume that $d < i_r(M)$. Then $d  \in \{2,4,\ldots,20\}$. We have 
\[
\frac{1}{4}d^2 + d  <  51 = n-r \quad \Longleftrightarrow \quad d \in \{2,4,6,8,10,12\}. 
\]
It follows from Proposition \ref{estimate} that  $i(M) \in \{14,16,18,20,22\}$. We proved in \cite{BOR19}, with different methods, that $i(M) = 22 = i_r(M)$.
\end{re}

\end {document}